\documentclass[a4paper,10pt]{amsart}
\usepackage{amsmath}
\usepackage{amssymb}
\usepackage{amsthm}
\usepackage[all]{xy}
\usepackage[latin1]{inputenc}        
\usepackage[dvips]{graphics}
\usepackage[dvips]{graphicx}
\usepackage{amscd}
\usepackage{mathrsfs}
\usepackage[mathcal]{eucal}
\usepackage{fullpage}
\author{Francesco Polizzi}
\title[On surfaces of general type]{On surfaces of general type with $p_g=q=1, \; K^2=3$}
\address{Dipartimento di Matematica, Universit{\`a} di Roma ``Tor Vergata'', Via della Ricerca Scientifica 1, 00133 Roma (Italy)}
\email{polizzi@.mat.uniroma2.it}
\date{\today}

\newtheorem{inizio}{Lemma}[section]
\newtheorem{theorem}[inizio]{Theorem}
\newtheorem{corollary}[inizio]{Corollary}
\newtheorem{proposition}[inizio]{Proposition}
\newtheorem{lemma}[inizio]{Lemma} 
 
\newtheorem{definition}[inizio]{Definition}
\newtheorem{remark}[inizio]{Remark}

\newtheorem{notation}[inizio]{Notation}
\newtheorem*{teo}{Theorem}
\newtheorem*{teoA}{Theorem A}
\newtheorem*{teoB}{Theorem B}
\newtheorem*{teoC}{Theorem C}

\newcommand{\hB}{\widehat{B}}

\newcommand{\F}{\mathbb{F}_2}

\newcommand{\hF}{\widehat{\mathbb{F}}_2}
\newcommand{\hS}{\widehat{S}}
\newcommand{\hL}{\widehat{L}}
\newcommand{\hmL}{\widehat{\mathcal{L}}}
\newcommand{\f}{\hat{f}_3}
\newcommand{\E}{\mathcal{E}}
\newcommand{\La}{\widehat{\Lambda}}
\newcommand{\Th}{\widehat{\Theta}}
\newcommand{\hmu}{\hat{\mu}}
\newcommand{\gD}{\mathfrak{D}_0}

\newcommand{\mO}{\mathcal{O}}
\newcommand{\hR}{\widehat{R}}
\newcommand{\hW}{\widehat{W}}

\newcommand{\Co}{C_{\xi_0}}
\newcommand{\Cone}{C_{\xi_1}}
\newcommand{\Ctwo}{C_{\xi_2}}
\newcommand{\Cthree}{C_{\xi_3}}

\begin{document}
 
\abstract
The moduli space $\mathscr{M}$ of surfaces of general type with $p_g=q=1, \; K^2=g=3$ (where $g$ is the genus of the Albanese fibration) was constructed by Catanese and Ciliberto in \cite{CaCi93}. In this paper we characterize the subvariety $\mathscr{M}_2 \subset \mathscr{M}$ corresponding to surfaces containing a genus $2$ pencil, and moreover we show that there exists a non-empty, dense subset $\mathscr{M}^0 \subset \mathscr{M}$ which parametrizes isomorphism classes of surfaces with birational bicanonical map.
\endabstract
\maketitle
\tableofcontents

\section{Introduction}
As the title suggests, this paper deals with surfaces of general type having both geometric genus and irregularity equal to $1$. The problem of classification of surfaces of general type is of exponential computational complexity; see \cite{Ca92}, \cite{Man97}, \cite{Ch96}; nevertheless, one can hope to classify at least those with small numerical invariants. It is well known that there exist surfaces of general type with $p_g=q=0$, the first example having been given by Godeaux in \cite{Go31}; later on, many other examples have been discovered and carefully investigated. On the other hand, since any surface of general type $S$ satisfies $\chi(\mathcal{O}_S) >0$, it follows that $q >0$ implies $p_g >0$; therefore surfaces of general type with $p_g=q=1$ are the irregular ones with the lowest geometric genus, and this explains their importance in the framework of classification. However, these surfaces are still mysterious, and only a few families have been hitherto discovered. If $S$ is irregular, then it satisfies Debarre's inequality $K^2 \geq 2p_g$, and a well-known result of Miyaoka and Yau implies $K^2 \leq 9 \chi(\mathcal{O}_S)$; therefore if $p_g=q=1$ one has $2 \leq K^2 \leq 9$. A complete classification has been obtained only in in the cases $K^2=2,3$: see \cite{Ca81}, \cite{CaCi91}, \cite{CaCi93}; for higher values of $K^2$ one merely knows some ``sporadic'' examples: see \cite{Ca99}, \cite{Pol03}. It is worth remarking that it is still unknown whether it is possible to construct surfaces of general type with $p_g=q=1, \; K^2=6,\; 7,\; 9$. \\
Let $S$ is a minimal surface of general type with $p_g=q=1$; since $q=1$, the Albanese map of $S$ is a connected fibration $\alpha \colon S \longrightarrow E$, where $E$ is an elliptic curve. We denote by $F$ a general fibre of $\alpha$, and by $g:=g(F)$ the genus of $F$. Already at the end of the XIX century, Castelnuovo and Enriques were aware of the feature that distinguishes the irregular surfaces from the regular ones, namely the existence of continuous systems of curves which are not linearly equivalent: see in particular Chapter IX of \cite{En49}. One can easily construct such a object in $S$ by considering the twists of the canonical bundle by the pullback via $\alpha$ of a degree zero line bundle over $E$. This means fixing a zero point $0 \in E$ and considering for any $t \in E$ the line bundle $K+t:= K+F_t-F_0$, where $F_t$ is the fibre of $\alpha$ over $t$. By Riemann Roch and semi-continuity one has $h^0(K+t)=1$ for general $t$, hence denoting by $C_t$ the unique curve belonging to the linear system $|K+t|$ one obtains a continuous, $1-$dimensional algebraic system of curves $\{C_t\}_{t \in E}$ parametrized by the elliptic curve $E$. It is called \emph{paracanonical system} and denoted by $\{K \}$; see Section \ref{surfaces} for further details. The paracanonical system is a powerful tool in the study of surfaces of general type with $p_g=q=1$, because it yields a natural map into a symmetric product of  $E$ in the following way. Let $\iota=\iota(K)$ be the \emph{index} of $\{ K \}$, namely the number of paracanonical curves passing through a general point of $S$; therefore we define the so called \emph{paracanonical map}
\begin{displaymath}             
\begin{split}  
&\omega \colon S \longrightarrow E(\iota)\\ 
&\omega(x):=t_1+ \cdots +t_{\iota}
\end{split}
\end{displaymath}
where $x \in S$ is a general point and $C_{t_1}, \ldots, C_{t_\iota}$ are the paracanonical curves containing it.\\
The symmetric product $E(n)$ is a $\mathbb{P}^{n-1}-$bundle over $E$, which is the projectification of the Atiyah bundle $\mathcal{E}_n$ of rank $n$ (see Section \ref{symmetric products}); let us denote by $F$ and $D$ the algebraic classes of a fiber and of a tautological section respectively. In the sequel $S$ always denotes a minimal surface of general type with $p_g=q=1, \; K^2=g=3$. 
The starting point for our investigations is (see Theorem \ref{centrale})
\begin{teo}[Catanese-Ciliberto]
If $E$ is an elliptic curve and $\mathfrak{D}$ is an effective divisor in
$E(3)$ algebraically equivalent to $4D-F$, then the general element of $|\mathfrak{D}|$ is a smooth surface with $p_g=q=1, \; K^2=g=3$. Conversely, given a smooth surface of general type with $p_g=q=1, \; K^2=g=3$, its paracanonical map $\omega \colon S \longrightarrow E(3)$ is a birational morphism onto its image $X:=\omega(S)$, which is isomorphic to the canonical model of $S$.  
The moduli space $\mathscr{M}$ of such surfaces is smooth, irreducible of dimension $5$.
\end{teo}
Up to translations, we can suppose $\mathfrak{D}=\mathfrak{D}_0$, where $\mathfrak{D}_0$ is linearly equivalent to $4D_0-F_0$; in this way the study of surfaces with $p_g=q=1, \; K^2=g=3$ is translated into the study of divisors in $|\mathfrak{D}_0|$ having at worst rational double points (R.D.P.) as singularities.\\
Let us briefly describe now the three parts in which the paper is divided. \\ \\
In Part 1 we collect some known results about the $n-$th symmetric product $E(n)$ of an elliptic curve $E$, and we use the Schr{\"o}dinger representation of $\mathcal{H}_3$,  the Heisenberg group of level $3$, described in \cite{CaCi93} and \cite{ADHPR93} in order to partially answer the following question that Ciliberto and Catanese pose in [CaCi93, Remark 3.4]: is the linear system $|\mathfrak{D}_0|$, which has at most simple base points, actually base-point free? Although we could not solve completely this problem, nevertheless we have been able to obtain the result below (cf. Theorem \ref{4basepoints}).
\begin{teoA}
The linear system $|\mathfrak{D}_0|$ contains at most the four base points:
\begin{displaymath}
0+\eta_1+ 2\eta_1, \ldots, 0+\eta_4+2\eta_4
\end{displaymath}
where $\eta_i$ is a 3-torsion point of $E$ different from $0$.
\end{teoA}
Part $2$ is devoted to the study of surfaces with $p_g=q=1, \; K^2=g=3$ and a pencil $|G|$ of curves of genus $2$. Xiao constructs them in \cite{Xi85b} by using Jacobian fibrations and modular curves, and he proves that their moduli space $\mathscr{M}_2$ is $1-$dimensional and it is isomorphic to the coarse moduli space of elliptic curves. Therefore one wonders whether it is possible to characterize $\mathscr{M}_2$ as a subvariety of the moduli space $\mathscr{M}$ described by Catanese and Ciliberto. This is equivalent to characterize the elements in $|\mathfrak{D}_0|$ with at worst R.D.P. as singularities and containing a genus $2$ pencil. Our solution is (cf. Theorem \ref{theorem 1})
\begin{teoB}
There is exactly one divisor $\mathcal{X} \in |\mathfrak{D}_0|$ containing a genus $2$ pencil, and it is smooth. Moreover, let $\xi_i$ be any $2-$torsion point of $E$ different from zero, and consider the curve:
\begin{displaymath} 
A_{\xi_i}= \{\xi_i+x+(\div x) \; | \;x \in E \}.
\end{displaymath}
Then we have $h^0(E(3), \; \mathfrak{D}_0 \otimes \mathscr{I}_{A_{\xi_i}})=1$, and $\mathcal{X}$ coincides with the only effective divisor in the linear system $|\mathfrak{D}_0 \otimes \mathscr{I}_{A_{\xi_i}}|$. 
\end{teoB}
The proof of Theorem B is based on a careful analysis of the reducible paracanonical divisors of $S$, in the case when $S$ contains a genus $2$ pencil: it turns out in particular that there is a very fine interplay between reducible paracanonical divisors and reducible elements of $|G|$: this is explained in Proposition \ref{decomposizione}, Corollary \ref{emmei} and Proposition \ref{riducibili}. Our method gives also some other interesting results: for instance, we are able to compute the intersection number $FG$, obtaining $FG=4$ (Proposition \ref{FG}). Finally, we give a description ``{\`a} la Campedelli'' of surfaces $S$ with $p_g=q=1, \; K^2=g=3$ containing a genus $2$ pencil, showing that they are double coverings of a Hirzebruch surface $\mathbb{F}_2$, branched over a singular curve that we describe in detail (Proposition \ref{covering}). \\ \\
In Part 3 we study the bicanonical map of surfaces of general type with $p_g=q=1, \; K^2=g=3$. Over the years, the problem of understanding the pluricanonical maps of surfaces of general type has attracted a considerable amount of attention; one of the most interesting cases is that of the bicanonical map $\phi:=\phi_{2K}$. The so-called \emph{standard case} for the non-birationality of the bicanonical map happens when one considers a surface $S$ with a genus $2$ pencil $|G|$: in this case, indeed, the bicanonical system of $S$ cuts out on the general element $G$ of the pencil a $g^2_4$ composed with the bicanonical involution, hence $\phi$ cannot be birational.    
There are only finitely many families of surfaces of general type whose bicanonical map is not birational and which present the non-standard case, i.e. such that they are not fibered by curves of genus $2$; unfortunately, their complete classification is still missing, see \cite{Ci97} for an account of this subject.
The behaviour of the bicanonical map for surfaces with $q(S) \geq 2$ or $\chi(\mathcal{O}_S) >1$ is completely understood by work of Catanese, Ciliberto, Francia, Mendes-Lopes, Borrelli: see \cite{CFM97}, \cite{CCM98}, \cite{CM02}, \cite{CM02b}, \cite{Bor03}. The remaining cases $p_g=q=0$ and $p_g=q=1$ are therefore the focus of current research in this area. The cases $p_g=q=0, \; K^2 \geq 7$ and $p_g=q=1, \; K^2=9$ have been settled out respectively by Pardini and Mendes Lopes in \cite{MP01} and \cite{MP03} and by Ciliberto and Mendes Lopes in \cite{CM02}, whereas the case $p_g=q=1, \; K^2=8$ has been settled out by the author in \cite{Pol03}. Here we deal with the case $p_g=q=1, \; K^2=g=3$, and our result is (cf. Theorem \ref{main3})
\begin{teoC} 
There exists a non-empty dense subset $\mathscr{M}^0 \subset \mathscr{M}$ which parametrizes isomorphism classes of surfaces with \emph{birational} bicanonical map. 
\end{teoC}
Singular paracanonical divisors of $S$ play a crucial role also in the proof of Theorem C: indeed we show that the bicanonical map of a smooth surface in $|\mathfrak{D}_0|$ is birational whenever the three paracanonical curves $C_{\xi}$, where $\xi$ is a non-zero $2-$torsion point of $E$, are smooth (Theorem \ref{main2}); this turns out to be a nonempty, open condition (Corollary \ref{insieme denso}). \\ \\
Finally, we want to point out the main problem, about the bicanonical map, left open by our results: namely the existence of minimal surfaces of general type with $p_g=q=1, \; K^2=3$ and presenting the non-standard case. In other words, does there exist some elements in $|\mathfrak{D}_0|$, different from $\mathcal{X}$, with at most R.D.P. as singularities and bicanonical map of degree $d>1$? We don't know any examples, albeit we have some results in this direction. In fact we showed that: 
\begin{itemize}
\item[-] the case $d=3$ never occurs (Proposition \ref{non3}). It is worth saying that the proof involves an old result of Kummer about quartic surfaces in $\mathbb{P}^3$ containing infinitely many conics (Theorem \ref{kummer});
\item[-] in the case $d=2$, if a surface $S$ presenting the non-standard case exists then the quotient of $S$ with respect to the bicanonical involution must be an elliptic surface with $p_g=1, \; q=0$ and Kodaira dimension $1$ (Proposition \ref{due casi grado 2}). This improves an earlier result of Xiao (Theorem \ref{Xiao list}).         
\end{itemize}
Another interesting problem is to describe the bicanonical image $\Sigma$ of a general surface with $p_g=q=1, \; K^2=g=3$; Theorem C tells us that $\Sigma$ has degree $12$ in $\mathbb{P}^3$. We solved the corresponding problem when $S$ contains a genus $2$ pencil: in this case the bicanonical map has degree $2$ and the bicanonical image is a surface of degree $6$ with a quadruple line and seven isolated nodes, see Proposition \ref{imag bic}. \\ \\  
$\mathbf{Notations \; and \; conventions.}$ All varieties, morphisms, etc. in this article are defined over the field $\mathbb{C}$ of the complex numbers. By ``surface'' we mean a projective, non-singular surface $S$, and for such a surface $K$ denotes the canonical class, $p_g(S)=h^0(S, \; K)$ is the \emph{geometric genus},  $q(S)=h^1(S, \; K)$ is the \emph{irregularity} and $\chi(\mathcal{O}_S)=1-q(S)+p_g(S)$ is the \emph{Euler characteristic}. We call $\phi$ the \emph{bicanonical map} of $S$, namely the rational map $\phi \colon S \longrightarrow \mathbb{P}^{K^2 +\chi(\mathcal{O}_S)-1}$ induced by the complete linear system $|2K_S|$. We denote by $X$ the canonical model of $S$, and we call again $\phi$ the bicanonical map of $X$. When the canonical class of $S$ is ample, that is when $S$ does not contains $(-2)-$curves, we often identify $S$ with $X$. If $S$ is a surface with $p_g=q=1$, then $\alpha \colon S \longrightarrow E$ is the Albanese map of $S$, and $g:=g(F)$ the genus of its general fibre $F$.
$\xi_1, \; \xi_2 , \; \xi_3$ are the non-zero $2-$torsion points of $E$, whereas $\eta_1, \cdots, \eta_8$ are the non-zero $3-$torsion points. The zero point in the law group of $E$ is usually denoted by $0$, but sometimes it is convenient for us to call it $\xi_0$ or $\eta_0$. If $X, \; Y$ are two varieties, then $X \cong Y$ means that $X$ and $Y$ are isomorphic. If $D_1, \; D_2$ are two divisors, then $D_1 \cong D_2$ means that $D_1$ and $D_2$ are linearly equivalent, whereas $D_1 \equiv D_2$ means that they are algebraically or numerically equivalent. \\ \\
$\mathbf{Acknowledgements.}$ This work is part of my Ph.D. Thesis at the University of Rome ``Tor Vergata''. I thank most sincerely Professor Ciro Ciliberto for his guidance and precious suggestions; without his help this paper would not have been possible. Part of this research was done during the spring semester of the academic year 2001-2002, when I was a visiting student at the Centro de Matem{\'a}tica e Aplica{\c c}{\~o}es Fundamentais (Universidade de Lisboa, Portugal), and I was supported by Eu Research Training Network EAGER, no. HPRN-CT-2000-00099. In particular I am indebted with M. Mendes Lopes for many useful discussions and for her supportive attitude. 

\part{Symmetric product of elliptic curves and surfaces of general type with $p_g=q=1$} \label{part:background}
\section{Symmetric products of elliptic curves} \label{symmetric products}
For completeness and clarity, we include in this section some basic facts about symmetric products of elliptic curves that we use in the whole paper; for statements which are given without proofs, we refer the reader to [CaCi93].\\
Let $E$ be an elliptic curve, and fix a zero point $0$ over it;
the $n$-th symmetric product $E(n)$ is the quotient of the
cartesian product $E^n$ by the natural action of the symmetric
group $S_n$; hence we can identify the elements of
$E(n)$ with the effective divisors of degree $n$ on $E$. We denote
by $\oplus$ the sum in the law group of $E$ and by $\div u$ the
opposite of the element $u \in E$. For any given $u$ we can
consider the divisors:
\begin{displaymath}
\begin{array}{cc}
D_u= \{u+x_2+ \cdots +x_n \thinspace | \thinspace x_2, \ldots, x_n \in E \};\\
F_u= \{x_1+ \cdots +x_n \thinspace | \thinspace x_1 \oplus \cdots
\oplus x_n=u \}.
\end{array}
\end{displaymath}
The algebraic equivalence classes of $D_u$ and $F_u$ are independent from $u$, and we denote them by $D$ and $F$ respectively. It is possible to show that these classes generate the N{\'e}ron-Severi group of $E(n)$, that is every divisor in $E(n)$ is algebraically equivalent to $aD+bF$ for some $a,b \in \mathbb{Z}$.\\ 
We observe moreover that the map sending $u+x_2+ \cdots +x_n$ to
$x_2+ \cdots x_n$ determines an isomorphism between $D_u$ and
$E(n-1)$; we will call it the \emph{natural identification} between
these two varieties.
\begin{inizio} \label{restringere}
Let us identify $D_u$ to $E(n-1)$ in the natural way. Therefore
we have:
\begin{enumerate}
\item $\mathcal{O}_{E(n)}(D_v) \otimes \mathcal{O}_{D_u}=
\mathcal{O}_{E(n-1)}(D_v)$; \item $\mathcal{O}_{E(n)}(F_v) \otimes
\mathcal{O}_{D_u}= \mathcal{O}_{E(n-1)}(F_{v \div u})$.
\end{enumerate}
\end{inizio}
\begin{proof}
Immediate.
\end{proof}
On the other hand $F_u$ is the fibre over $u$ of the Abel-Jacobi
map of $E(n)$:
\begin{displaymath}
\begin{array}{cc}
\beta \colon E(n) \longrightarrow E\\
\beta (u_1+ \cdots +u_n):=u_1 \oplus \cdots \oplus u_n.
\end{array}
\end{displaymath}
Therefore $F_u$ is isomorphic to $\mathbb{P}^{n-1}$ ([ACGH85, p.18]), that is $E(n)$ is a $\mathbb{P}^{n-1}-$bundle over $E$.
Since $H^2(E, \mathcal{O}_E)=0$, it follows that $E(n)$ is the
projectification $\mathbb{P}(\mathcal{E}_n)$ of a suitable rank $n$
vector bundle $\mathcal{E}_n$ over $E$. The vector bundle
$\mathcal{E}_n$ can be constructed in the following inductive way:
$\mathcal{E}_1=\mathcal{O}_E(0)$, $\mathcal{E}_2$ is defined as 
the only nontrivial extension of line bundles:
\begin{displaymath}
0 \longrightarrow \mathcal{O}_E \longrightarrow \mathcal{E}_2 \longrightarrow
\mathcal{O}_E(0) \longrightarrow 0
\end{displaymath}
and $\mathcal{E}_n$ is defined as the only nontrivial extension of vector bundles:
\begin{displaymath}
0 \longrightarrow \mathcal{O}_E \longrightarrow \mathcal{E}_n \longrightarrow
\mathcal{E}_{n-1} \longrightarrow 0.
\end{displaymath}
We call $\mathcal{E}_n$ the \emph{Atiyah bundle} of rank $n$.
From the previous construction it follows that $K_{E(n)} \cong -nD_0+F_0$.\\
In particular $E(2)$ is a ruled surface over the elliptic curve $E$ with invariant $e=-1$ (for the definition of $e$, see [Ha77, Chapter V]). In this case, sometimes we will denote by $f$ and $h$ the class of the fibre and of the section respectively; this is done for the first time in Remark \ref{altropuntodivista}.\\ 
We observe that any automorphism $g$ of $E$ naturally induces an
automorphism $g(n)$ of $E(n)$; hence $E(n)$ has a
$1-$parameter group of automorphisms induced by the group of
translations of $E$. If $u \in E$, we see that $u$ acts on the set
of classes of algebraically equivalent divisors by $D_t
\longrightarrow D_{t \oplus u}$, $F_t \longrightarrow F_{t \oplus nu}$.
\begin{proposition} \label{classi}
The group of the translations of $E$ acts transitively on all algebraically equivalent classes of linear equivalence which are not an integral multiple of the canonical class.\\
Moreover the homomorphism $g \longrightarrow g(n)$ is an isomorphism
of \textrm{Aut}$(E)$ onto \textrm{Aut}$(E(n))$.
\end{proposition}
\begin{proof}
See [CaCi93, Proposition 1.5].
\end{proof}
An immediate consequence of Proposition \ref{classi} is that the
cohomology of a divisor which is not algebraically equivalent to
some multiple of $-nD+F$ depends only on its algebraic equivalence
class; indeed, we have the following:
\begin{theorem}[Catanese, Ciliberto] \label{coomologia in E(n)}
Let $\mathfrak{D}$ be a line bundle associated with a divisor
algebraically equivalent to $aD+bF$. Then:
\begin{enumerate}
\item if $a+nb \neq 0$, there is exactly one nonzero cohomology
group of $\mathfrak{D}$, in the following cases:
$H^0(\mathfrak{D})$ if $a \geq 0$, $a+nb >0$, $H^1(\mathfrak{D})$
if $a \geq 0, \; a+nb <0$, $H^{n-1}(\mathfrak{D})$ if $a \leq -n,
 \;a+nb >0$, $H^n(\mathfrak{D})$ if $a \leq -n, \; a+nb <0$; \item if
$-n <a <0$ all cohomology groups vanish;
\item if $a+nb=0$ then the Euler-Poincar{\'e} characteristic of $\mathfrak{D}$ is zero and there are at most two nonzero cohomology groups, namely\\\\
$H^0(\mathfrak{D})$ and $H^1(\mathfrak{D})$ if $a \geq 0$, \quad $H^{n-1}(\mathfrak{D})$ and $H^n(\mathfrak{D})$ if $a \leq -n$ \\\\
and only for finitely many divisors $\mathfrak{D}$ algebraically equivalent to $aD+bF$ are these cohomology groups nonzero.
\end{enumerate}
Moreover, the Euler-Poincar{\'e} characteristic of $\mathfrak{D}$ is:
\begin{displaymath}
\chi(\mathfrak{D})=(n!)^{-1}(a+nb)\cdot \prod_{i=1}^{n-1}(a+i).
\end{displaymath}
\end{theorem}
\begin{proof}
See [CaCi93, Theorem  1.17].
\end{proof}
The cohomology of the divisors algebraically equivalent to a
multiple of the canonical class has a particular behavior,
since it depends also on the linear equivalence class. By using the theory of representation of the Heisenberg group of level $n$, it is
possible to prove the following:
\begin{proposition}[cf. \cite{CaCi93}, Section 2] \label{plurianticanonici}
The only effective classes in $E(n)$ which are algebraically equivalent to a multiple of the canonical class are those obtained by adding to a
plurianticanonical class a divisor class of $n-$torsion.
\end{proposition}
Here we give some examples:
\begin{itemize}
\item $n=2$.\\
We have $h^0(E(2), -K_{E(2)})=h^0(E(2), \;  2D_0-F_0)=0$, whereas $h^0(E(2), -K_{E(2)}+\xi)=1$ if $\xi$ is a $2-$torsion point different from $0$. 
Hence we obtain three effective divisors algebraically but not linearly equivalent to the anticanonical divisor of $E(2)$; they are the three smooth elliptic curve given by:
\begin{displaymath}
T_{\xi_i}:= \{t+(t \oplus \xi_i) \thinspace |\thinspace t \in E \}
\cong -K_{E(2)}+\xi_i,
\end{displaymath}
where $\xi_i, \; i=1,2,3$ are the three non zero $3-$torsion points of $E$.
In the same way, $h^0(E(2), -2K_{E(2)})=2$  but $h^0(E(2),
-2K_{E(2)}+\xi_i)=1$. The elements of the anti-bicanonical pencil $|-2K_{E(2)}|$ are the curves $T_a:=\{x+(x\oplus a) \; | \; x \in E \;  \}$ with $2a \neq 0$; notice that $T_a=T_{\div a}$. Again, we have three curves algebraically but
not linearly equivalent to the divisor $ -2K_{E(2)}$, namely the
three reducible curves $T_{\xi_i}+T_{\xi_j}$, where $\xi_i$ and $\xi_j$
are two distinct nonzero points of $2-$torsion of $E$. In fact, one has:
\begin{displaymath}
T_{\xi_i}+T_{\xi_j} \cong -2K_{E(2)}+(\xi_i \oplus \xi_j).
\end{displaymath}
\item[]
\item $n=3$.\\
We have $h^0(E(3), -K_{E(3)})=2$; the anticanonical pencil is
generated by two of the three anticanonical divisors $\Psi_{\xi}$,
where $\xi$ is a nonzero $2-$torsion points of $E$ and:
$$\Psi_{\xi}:=\{u+v+(u \oplus \xi) \thinspace | \thinspace u,v \in E \}.$$
The base locus of $|-K_{E(3)}|$ is the curve:
$$\Gamma:= \{x+(x \oplus \xi_i)+(x \oplus \xi_j) \thinspace | \thinspace x \in E\},$$
where $\xi_i$, $\xi_j$ are two distinct nonzero $2-$torsion points of $E$. Notice that $\Gamma$ does not depends on the choice of $\xi_i$ and $\xi_j$.\\
Moreover, we have $h^0(E(3), -K_{E(3)}+ \eta)=1$, where $\eta$ is
one of the eight $3-$torsion point in $E-\{0 \}$. This means that
$E(3)$ contains exactly $8$ divisors algebraically but not
linearly equivalent to $-K_{E(3)}$, corresponding to the $0-$dimensional linear systems $|-K_{E(3)}+\eta_i|, \; i=1, \ldots,8$. We describe these divisors in Section \ref{basepointsD0}.
\end{itemize}
\section{Surfaces of general type with $p_g=q=1$} \label{surfaces}
Let $S$ be a minimal surface of general type. If $p_g=0$, then by using the inequality $\chi(\mO_S)>0$ we obtain $q=0$; this shows that the surfaces of general type with $p_g=q=1$ are the irregular ones with the lowest geometric genus. It is well known that for any irregular surface of general type Debarre's inequality $K^2 \geq 2p_g$ holds (see [De82]); on the other hand, Miyaoka's inequality yields $K^2 \leq 9 \chi(\mO_S)$ for any minimal surface of general type, then if $p_g=q=1$ we have $2 \leq K^2 \leq 9$. Surfaces with $p_g=q=1, \; K^2=2$ were classified in [Ca81], whereas [CaCi91] and [CaCi93] deal with the case $K^2=3$. For higher values of $K^2$ only some ``sporadic'' examples are known: see [Ca99] and [Pol03]. \\
When $p_g=q=1$, there are two basic tools that one can use in order to study the geometry of $S$: the \emph{Albanese fibration} and the \emph{paracanonical system}.  First of all, $q=1$ implies that the Albanese variety of $S$ is an elliptic curve $E$, hence the Albanese map $\alpha \colon S \longrightarrow E$ is a connected fibration; we denote by $F$ the general fibre of $\alpha$ and by $g=g(F)$ its genus. Let us fix a zero point $0 \in E$, and for any $t \in E$ let us write $K+t$ for the line bundle $K+F_t-F_0$. By Riemann-Roch we obtain: 
\begin{displaymath}
h^0(S,\; K+t)=1+h^1(S, \; K+t) 
\end{displaymath}
for any $t \in E- \{ 0 \}$. Since $p_g=1$, by semicontinuity there is a Zariski open set $E' \subset E$, containing $0$, such that for any $t \in E'$ we have $h^0(S, \;K+t)=1$; we denote by $C_t$ the unique curve in $|K+t|$. The \emph{paracanonical incidence correspondence} is defined to be the surface $Y$ in $S \times E$ which is the scheme-theoretical closure of the set $ \{(x,t) \in S \times E' \; | \; x \in C_t \}$. Then one can define $C_t$ for any $t$ as the fibre of $Y \longrightarrow E$ over $t$, and $Y$ provides a flat family of curves on $S$, that we denote by $\{K \}$ or by $\{C_t \}$ and we call the \emph{paracanonical system} of $S$. According to [Be88], $\{K \}$ is the irreducible component of the Hilbert scheme of curves on $S$ algebraically equivalent to $K$ which dominates $E$. We write $\{K \}=Z + \{M \}$, where $Z$ is the \emph{fixed} part of $\{K \}$ and $\{ M \}$ is the \emph{movable} part. Let $\mathcal{P}$ be a Poincar{\'e} sheaf on $S \times E$; then we call $\mathcal{K}=\pi_S^*(\omega_S) \otimes \mathcal{P}$ the paracanonical system on $S \times E$. Let $\Lambda^i:= R^i (\pi_E)_* \mathcal{K}$. By the base change theorem, $\Lambda^0$ is an invertible sheaf on $E$,  $\Lambda^2$ is a skyscraper sheaf of length $1$ supported at the origin and $\Lambda^1$ is zero at the origin, and supported on the set of points $\{ t \in E \; | \; h^0(S, \;K+t) > 1 \}$; let $\lambda:= \textrm{length}(\Lambda^1)$. 
\begin{definition}
The \emph{index} $\iota=\iota(K)$ of the paracanonical system is the intersection number $Y \cdot (\{ x \} \times E)$. Roughly speaking, $\iota$ is the number of paracanonical curves through a general point of $S$.
\end{definition}
If $g=g(F)$ is the genus of the Albanese fibre of $S$, then the following relation holds:
\begin{displaymath} 
\iota=g - \lambda.     
\end{displaymath}
$\lambda$ is actually a topological invariant of $S$, see [CaCi91, Theorem 1.4]. Moreover $\iota=1$ if and only if the movable part $\{ M \}$ of the paracanonical system coincides with the Albanese fibration $\{ F \}$. 
 
The \emph{paracanonical map} 
\begin{displaymath}
\omega \colon S \longrightarrow E(\iota) 
\end{displaymath}
is defined in the following way: if $x$ is a general point of $S$, then $\omega(x):=t_1+ \cdots +t_{\iota}$, where $C_{t_1}, \ldots, C_{t_{\iota}}$ are the paracanonical curves of $S$ passing through $x$. 
There is a commutative diagram:
\begin{displaymath}
\xymatrix{
{S} \ar[dr]_{\omega} \ar[r]^{\alpha} 
& {E} \\
& E(\iota)  \ar[u]_{\beta} }
\end{displaymath}
When $\omega$ is birational, the fibres of $\beta$ cut out on $\omega(S)$ the Albanese fibration, whereas the divisors of $\{ D \}$ cut out the paracanonical system. \\
The best result that one might obtain would be to classify the triples $(K^2,g, \iota)$ such that there exists a minimal surface of general type $S$ with $p_g=q=1$ and these invariants. Since by the results of Gieseker the moduli space $\mathscr{M}_{\chi, \; K^2}$ of surfaces of general type with fixed $\chi(\mathcal{O}_S), \; K^2$ is a quasiprojective variety, it turns out that there exist only finitely many such triples, but a complete classification is still missing. We already said that only the cases $K^2=2, \; 3$ are hitherto completely understood and what is known is:
\begin{theorem}[Catanese] \label{cataneseK2}
Let $S$ be a minimal surface of general type with $p_g=q=1, \; K^2=2$. Then one has $g=\iota=2$ and the paracanonical map, coinciding with the relative canonical map, is a $2:1$ morphism $\omega \colon S \longrightarrow E(2)$, branched along a divisor algebraically equivalent to $6D-2F$ $($see Section \ref{symmetric products} for a description of the N{\'e}ron-Severi group of $E(n)$$)$. The general paracanonical curve of $S$ is irreducible. 
\end{theorem} 
\begin{proof}
See [Ca81]
\end{proof}
\begin{theorem}[Catanese, Ciliberto] \label{catanese-ciliberto}
Let $S$ be a minimal surface of general type with $p_g=q=1, \; K^2=3$. Then either $g= \iota=2$ or $g= \iota=3$, and in both cases the paracanonical map coincides with the relative paracanonical map.
More precisely:
\begin{itemize}
\item in the case $g=\iota=2$, $\omega \colon S \longrightarrow E(2)$ is not a morphism but its indeterminacy can be resolved by blowing up a single point of $S$. Moreover $\omega$ is generically $2:1$ and the branch curve is algebraically equivalent to $6D$. 
\item in the case $g=\iota=3$, $\omega  \colon S  \longrightarrow E(3)$ is a morphism which is birational onto its image. Moreover, $\omega$ is an isomorphism of the canonical model of $S$ onto $\omega(S)$, which is a divisor with at most simple singularities in a linear system $|\mathfrak{D}|$, where $\mathfrak{D}$ is a divisor homologous to $4D-F$.
\end{itemize} 
In any case the general paracanonical curve of $S$ is irreducible.
\end{theorem}
\begin{proof}
See [CaCi91], [CaCi93].
\end{proof}
Let us consider now in greater detail the case $K^2=g=3$. In this situation one has $\iota=3$, hence $\lambda=0$; this means $h^0(K+t)=1$ for any $t \in E$. Theorem \ref{catanese-ciliberto} says that the canonical model of $S$, that we call $X$, is a divisor in $E(3)$ algebraically equivalent to $4D-F$; more precisely, there is the following result that we can consider the starting point for our investigations.
\begin{theorem}[CaCi93, Theorem 3.1] \label{centrale}
If $E$ is an elliptic curve and $\mathfrak{D}$ is an effective divisor in
$E(3)$ algebraically equivalent to $4D-F$, then the general element of $|\mathfrak{D}|$ is a
smooth surface with $p_g=q=1, \; K^2=g=3$. Conversely, if we have a smooth surface of general type with $p_g=q=1, \; K^2=g=3$, then its paracanonical map $\omega \colon S \longrightarrow E(3)$ is a birational morphism onto its image $\omega(S)$, which is isomorphic to the canonical model $X$ of $S$.      
The corresponding moduli space $\mathscr{M}$ is generically smooth, irreducible, uniruled of dimension $5$.
\end{theorem} 
Notice that we have $|\mathfrak{D}|=|-K_{E(3)}+D_u|$ for some $u \in E$.
\begin{remark}
The dimension of $\mathscr{M}$ can be computed as follows. Since the group of translations of $E(3)$ acts transitively on the set of linear equivalence classes of divisor algebraically equivalent to $4D-F$ $($Proposition \ref{classi}$ \;)$, we can suppose that the canonical model $X$ of $S$ belongs to $|\mathfrak{D}_0|$, where $\mathfrak{D}_0$ is a divisor linearly equivalent to $4D_0-F_0$. Since $h^0(E(3), \; \mathfrak{D}_0)=5$ $($Theorem \ref{coomologia in E(n)} $)$ and the elliptic curves depend on $1$ parameter, we obtain $\textrm{dim} \;\mathscr{M}=5$.
\end{remark} 
 Vice versa, any element in $|\mathfrak{D}_0|$ with at most rational double points as singularities is the canonical model of a minimal surface of general type with $p_g=q=1, \; K^2=g=3$. Moreover, the linear system $|\mathfrak{D}_0|$ contains at most simple base points ([CaCi93, Lemma 3.3]); it follows by Bertini's theorem that the general element of $|\mathfrak{D}_0|$ is smooth. This in turn implies that the general surface $S$ with $p_g=q=1, \; K^2=g=3$ has ample canonical class. In this case we often identify it with its canonical model $X \in |\gD|$.
\section{Base points of $|\mathfrak{D}_0|$} \label{basepointsD0}
In [CaCi93, Remark 3.4] the authors wonder whether the linear system $|\mathfrak{D}_0|$, which has at most isolated, simple base points (cf. [CaCi93, Proposition 1.19 and Lemma 3.3]), is actually base-point free. They suspect so, but don't give a complete proof. In this section we prove that $|\mathfrak{D}_0|$ has at most four base points (see Theorem \ref{4basepoints}). This result was inspired by the paper \cite{ADHPR93}, which contains a careful analysis of the surfaces belonging to the linear systems $|-K_{E(3)}|$ and $|-K_{E(3)}+\eta_i|$. Let us denote by $\mathcal{X}_{\eta_i}$ the unique element of $|-K_{E(3)}+\eta_i|$, and let $\mathcal{A} \in |-K_{E(3)}|$.     
\begin{proposition}[ADHPR93]\label{ADHPR}
The following holds:
\begin{enumerate}
\item the general element $\mathcal{A}$ in the pencil $|-K_{E(3)}|$ is a smooth, abelian surface isogenous to $E \times E$, and the restriction of $\beta \colon E(3) \longrightarrow E$ to $\mathcal{A}$ gives rise to a smooth, isotrivial elliptic fibration $\mathcal{A} \longrightarrow E$. Moreover the pencil contains four singular elliptic scrolls $\mathcal{A}_1, \ldots, \mathcal{A}_4$;
\item any $\mathcal{X}_{\eta_i}$ is a smooth bielliptic surface of type $6$ $($see Chapter VI of \cite{Be96}$)$. The restriction of $\beta$ gives rise to a smooth, isotrivial elliptic fibration $\mathcal{X}_{\eta_i} \longrightarrow E$, whereas the pencil $|-K_{E(3)}|$ cuts out on $\mathcal{X}_{\eta_i}$ an isotrivial pencil of elliptic curves with three triple fibres.
\end{enumerate}   
\end{proposition}
The proof of this result, as well as the proof of Theorem \ref{coomologia in E(n)}, is based on the representation theory of the Heisenberg group. Here we are mostly interested in the Heisenberg group of level $3$, which is the central extension:
\begin{displaymath}   
1 \longrightarrow \mu_3 \longrightarrow \mathcal{H}_3 \longrightarrow \mathbb{Z}_3 \times \mathbb{Z}_3 \longrightarrow 1
\end{displaymath}
defined by the relation $[\sigma,\tau]= \varepsilon^{-1}$, where $\sigma, \; \tau$ generate $\mathbb{Z}_3 \times \mathbb{Z}_3$ and $\varepsilon:=e^{\frac{2 \pi i}{3}}$ generates $\mu_3$. If $x_0, \; x_1, \; x_2$ are a basis of the vector space $H^0(\mathcal{O}_{\mathbb{P}^2}(1))$, the \emph{Schr{\"o}dinger representation} of $\mathcal{H}_3$ is defined as follows:
\begin{displaymath}
\sigma(x_i)=x_{i-1}; \quad \tau(x_i)=\varepsilon^{-i}x_i; \quad \varepsilon(x_i)=\varepsilon x_i,
\end{displaymath}
where $i$ is counted modulo $3$. If we consider the induced action on 
$H^0(\mathcal{O}_{\mathbb{P}^2}(3))$, the center $\varepsilon$ of $\mathcal{H}_3$ acts trivially, hence it descends to a representation of $\mathbb{Z}_3 \times \mathbb{Z}_3$, which splits in $1-$dimensional subrepresentations because $\mathbb{Z}_3 \times \mathbb{Z}_3$ is abelian. Let $(a,b)$ be the character of $\mathbb{Z}_3 \times \mathbb{Z}_3$ such that $\sigma$ acts by $\varepsilon^a$ and $\tau$ acts by $\varepsilon^b$, and let us write:
\begin{displaymath}
H^0(\mathcal{O}_{\mathbb{P}^2}(3))=\bigoplus_{(a,b)}V_{(a,b)},
\end{displaymath}
where $V_{(a,b)}$ is the eigenspace corresponding to the character $(a,b)$. The invariant subspace $V_{(0,0)}$ has dimension $2$ and it is generated by the two polynomials:
\begin{displaymath}
x_0^3+x_1^3+x_2^3, \quad x_0x_1x_2; 
\end{displaymath}
they span a pencil of cubics, called the \emph{Hesse pencil}, having the property that its base points coincide with the inflexion points of its elements. The Hesse pencil contains four singular members, namely the four triangles:
\begin{displaymath}
\begin{split}
T_{(1,0)}&=x_0x_1x_2\\    
T_{(1,2)}&=(x_0+\varepsilon^2 x_1+\varepsilon^2 x_2)(x_0+x_1+\varepsilon x_2)(x_0+\varepsilon x_1+x_2)\\
T_{(0,1)}&=(x_0+\varepsilon x_1+\varepsilon^2 x_2)(x_0+ \varepsilon^2 x_1+\varepsilon x_2)(x_0+x_1+x_2)\\
 T_{(1,1)}&=(x_0+ x_1+\varepsilon^2 x_2)(x_0+ \varepsilon x_1+\varepsilon x_2)(x_0+\varepsilon^2 x_1+x_2);
\end{split}  
\end{displaymath}
notice that the twelve edges of these triangles are the twelve lines joining three of the base points. All the subspaces $V_{(a,b)}$ with $(a,b) \neq (0,0)$ are $1-$dimensional, and their generators are the following:
\begin{equation} \label{cubichecaratteri}
\begin{split}
C_{(1,0)}&=x_0^3+\varepsilon x_1^3+\varepsilon^2 x_2^3\\
C_{(2,0)}&=x_0^3+ \varepsilon^2 x_1^3 +\varepsilon x_2^3\\
C_{(0,1)}&=x_0x_1^2+x_1x_2^2+x_2x_0^2\\
C_{(1,1)}&=x_0x_1^2+ \varepsilon x_1x_2^2+ \varepsilon^2 x_2x_0^2\\
C_{(2,1)}&=x_0x_1^2+ \varepsilon^2 x_1x_2^2+ \varepsilon x_2x_0^2\\
C_{(0,2)}&=x_0^2x_1+x_1^2x_2+x_2^2x_0\\
C_{(1,2)}&=x_0^2x_1+ \varepsilon x_1^2x_2+ \varepsilon^2 x_2^2x_0\\
C_{(2,2)}&=x_0^2x_1+ \varepsilon^2 x_1^2x_2+ \varepsilon x_2^2x_0.
\end{split}
\end{equation}
Propositions \ref{inters cubiche1} and \ref{inters cubiche2} describe how the curves $C_{(a,b)}$, $T_{(i,j)}$ intersect each other; since their proofs are nothing but straightforward computations, we  leave them to the reader.  
\begin{proposition} \label{inters cubiche1}
The curve $C_{(a,b)}$ contains the vertices of the triangle $T_{(i,j)}$ if and only if $(i,j) \neq \pm (a,b)$. In this case, the intersection multiplicity of $C_{(a,b)}$ and $T_{(i,j)}$ at each vertex is equal to $3$. 
\end{proposition}
\begin{proposition} \label{inters cubiche2}
Let $C_{(a_1,b_1)}, \; C_{(a_2,b_2)}$ be two different cubics in $($\ref{cubichecaratteri}$)$. Then all their intersection points are vertices of some triangle $T_{(i,j)}$, according to one of the following two cases:
\begin{enumerate}
\item if $(a_1,b_1)=-(a_2,b_2)$ then $C_{(a_1,b_1)}$ and $C_{(a_2,b_2)}$ intersect transversally at the vertices of the three triangles $T_{(i,j)}$ with $(i,j) \neq \pm (a_1,b_1)$;
\item if $(a_1,b_1) \neq -(a_2,b_2)$, then $C_{(a_1,b_1)}$ and $C_{(a_2,b_2)}$ intersect at the vertices of the two triangles $T_{(i,j)}$ with $(i,j) \neq \pm (a_1,b_1), \pm (a_2,b_2)$. More precisely, they intersect transversally at the vertices of $T_{\pm(a_1+a_2,b_1+b_2)}$ and they have a simple contact at the vertices of $T_{\pm(a_1-a_2, b_1-b_2)}$.   
\end{enumerate}
\end{proposition}
Since $[\sigma,\tau]=\varepsilon^{-1}$ on $H^0(\mathcal{O}_{\mathbb{P}^2}(1))$, the action of $\mathcal{H}_3$ descends to an action of $\mathbb{Z}_3 \times \mathbb{Z}_3$ on $\mathbb{P}(H^0(\mathcal{O}_{\mathbb{P}^2}(1))^*) \cong \mathbb{P}^2$. It is possible to give an action of $\mathbb{Z}_3 \times \mathbb{Z}_3$ on $E(3)=\mathbb{P}(\mathcal{E}_3)$, in such a way that on each fiber of this projective bundle it lifts to the Schr{\"o}dinger representation of $\mathcal{H}_3$; this can be done by identifying $\mathbb{Z}_3 \times \mathbb{Z}_3$ with the group of $3-$torsion points of $E$, and by considering the corresponding subgroup of translations: that is, $\eta_i$ acts on $x+y+z \in E(3)$ as follows: 
\begin{equation} \label{azione traslazione}    
\eta_i \cdot (x+y+z):=(x \oplus\eta_i)+(y \oplus \eta_i)+(z\oplus \eta_i).
\end{equation}
 Notice that, for any $i\neq 0$, the subgroup $\langle \; 1,\eta_i, 2\eta_i \; \rangle$ generated by $\eta_i$ fixes pointwise the curve $N_{\eta_i}:=\{x+(x \oplus \eta_i)+ (x\oplus 2 \eta_i)\; | \; x \in E \}$.    
\begin{proposition} \label{enneeta}
For any $i \neq 0$, $N_{\eta_i}$ is a smooth, elliptic curve in $E(3)$ such that:
\begin{displaymath}
N_{\eta_i}D=1, \quad N_{\eta_i}F=3.
\end{displaymath}
\end{proposition}
\begin{proof}
Let $\eta:=\eta_i$. Since $N_{\eta}$ is isomorphic to $E$, it is a smooth, elliptic curve. Set theoretically we can write:
\begin{displaymath}
N_{\eta} \cap D_0=N_{\eta} \cap D_{\eta}=N_{\eta} \cap D_{2\eta}=\{0+\eta+2 \eta\},
\end{displaymath}
and since $D_0,\; D_{\eta}, \; D_{2\eta}$ intersect transversally at the point $0+\eta+2\eta$, this implies $N_{\eta}D=1$. Let us consider now the surface $\Psi_{\xi}:=\{u+v+(u \oplus \xi)\; | \; u,v \in E\}$, where $\xi$ is one of the three non-zero $2-$torsion points of $E$. We know by Section \ref{symmetric products} that $\Psi_{\xi}$ belongs to the anticanonical pencil of $E(3)$, and it is immediate to verify that $\Psi_{\xi} \cap N_{\eta}=\emptyset$; hence $(3D-F)N_{\eta}=0$, that is $N_{\eta}F=3$.  
\end{proof}
The elements $\mathcal{A}\in |-K_{E(3)}|$ and the eight divisors $\mathcal{X}_{\eta_i}$ are invariant for the action (\ref{azione traslazione}) (see also Proposition \ref{classi}); if we restrict them to any fibre $F \cong \mathbb{P}^2$ of $E(3)$, then $|-K_{E(3)}|$ cuts out the Hesse pencil, whereas the $\mathcal{X}_{\eta_i}$'s cut out the eight cubics $C_{(a,b)}$. The four triangles $T_{(i,j)}$ are cut out on $F$ by the four singular scrolls $\mathcal{A}_1, \ldots, \mathcal{A}_4$. Moreover, Proposition \ref{enneeta} tells us that, for any nonzero $3-$torsion point $\eta$, the curve $N_{\eta}$ intersects any fibre $F$ of $E(3)$ transversally in $3$ points, which in turn correspond to the vertices of one of the triangles $T_{(i,j)}$. The nine base points of the Hesse pencil are cut out on $F$ by the curve $\Gamma$, base locus of $|-K_{E(3)}|$.
Looking at Propositions \ref{inters cubiche1} and \ref{inters cubiche2}, we obtain:
\begin{proposition} \label{inters biellittiche1}   
Any bielliptic surface $\mathcal{X}_{\eta_i}$ contains exactly three curves $N_{\eta_j}$, namely the ones with $\eta_j \neq \eta_i, 2\eta_i$. These curves are the three triple fibres of the pencil cut out on $\mathcal{X}_{\eta_i}$ by $|-K_{E(3)}|$. More precisely, each curve $N_{\eta_j}$ is contained in exactly one singular elliptic scroll $\mathcal{A}_j$, and it is the set-theoretical intersection of $\mathcal{A}_j$ and $\mathcal{X}_{\eta_i}$.   
\end{proposition}
\begin{remark}
The curve $N_{\eta_j}$ actually coincides with the singular locus of the elliptic scroll $\mathcal{A}_j$. Moreover, the scheme-theoretical intersection of $\mathcal{A}_j$ and $\mathcal{X}_{\eta_i}$ is not reduced, because $N_{\eta_j}$ appears in it with multiplicity $3$.
\end{remark}
\begin{proposition} \label{inters biellittiche2}
Let $\eta_i,\eta_j$ be two distinct, non-zero $3-$torsion points of $E$. Then the two surfaces $\mathcal{X}_{\eta_i}, \;\mathcal{X}_{\eta_j}$ intersect along some of the curves $N_{\eta_k}$, according to one of the following two cases:
\begin{enumerate}
\item if $\eta_j=2\eta_i$, then $\mathcal{X}_{\eta_i}$ and $\mathcal{X}_{\eta_j}$ intersect transversally along the three curves $N_{\eta_k}$ with $\eta_k \neq \eta_i, 2\eta_i$;
\item if $\eta_j \neq 2\eta_i$, then $\mathcal{X}_{\eta_i}$ and $\mathcal{X}_{\eta_j}$ intersect along the two curves $N_{\eta_k}$ with  $\eta_k \neq \eta_i, 2\eta_i, \eta_j, 2\eta_j$. More precisely, they intersect transversally along the curve $N_{\eta_i \oplus \eta_j}$, and they have a simple contact along the curve $N_{2(\eta_i \oplus \eta_j)}$.
\end{enumerate}
\end{proposition}
In the sequel we identify the group of $3-$torsion points of $E$ with $\mathbb{Z}_3 \times \mathbb{Z}_3$  in the following way:
\begin{displaymath}
\eta_1 \longleftrightarrow (0,1), \; \eta_2 \longleftrightarrow (1,0), \; \eta_3 \longleftrightarrow (1,1), \; \eta_4 \longleftrightarrow (1,2).
\end{displaymath}
Therefore, according with Proposition \ref{inters biellittiche2}, any intersection of $\mathcal{X}_{\eta_i}$ and $\mathcal{X}_{\eta_j}$ is described in the following table:
\begin{equation} \label{tavola int biellittiche}
\begin{tabular}{|c|c|c|c|c|c|c|c|c|}
\hline
 {} & $\mathcal{X}_{\eta_1}$ & $\mathcal{X}_{\eta_2}$ & $\mathcal{X}_{\eta_3}$ &  $\mathcal{X}_{\eta_4}$  &  $\mathcal{X}_{2\eta_1}$  &  $\mathcal{X}_{2\eta_2}$ &  $\mathcal{X}_{2\eta_3}$ &  $\mathcal{X}_{2\eta_4} $  \\
\hline
 $\mathcal{X}_{\eta_1}$ & $\bullet$ & $\eta_3,2\eta_4$ & $\eta_4,2\eta_2$ & $\eta_2, 2\eta_3$ & $\eta_2, \eta_3, \eta_4$ & $\eta_4, 2\eta_3$ & $\eta_2, 2\eta_4$ & $\eta_3, 2\eta_2$ \\
\hline
 $\mathcal{X}_{\eta_2}$ & $\eta_3, 2\eta_4$ & $\bullet$ & $\eta_4, 2\eta_1$ & $\eta_3, 2\eta_1$ & $\eta_4, 2\eta_3$ & $\eta_1,\eta_3, \eta_4$ & $\eta_1, 2\eta_4$ & $\eta_1, 2\eta_3$ \\
\hline
$\mathcal{X}_{\eta_3}$ & $\eta_4, 2\eta_2$ & $\eta_4, 2\eta_1$ & $\bullet$ & $\eta_2, 2\eta_1$ & $\eta_2, 2\eta_4$ & $\eta_1, 2\eta_4$ & $\eta_1, \eta_2, \eta_4$ & $\eta_1, 2\eta_2$ \\
\hline
$\mathcal{X}_{\eta_4}$ & $\eta_2, 2\eta_3$ & $\eta_3, 2\eta_1$ & $\eta_2, 2\eta_1$ & $\bullet$ & $\eta_3, 2\eta_2$ & $\eta_1, 2\eta_3$ & $\eta_1, 2\eta_2$ & $\eta_1, \eta_2, \eta_3$ \\
\hline
 $\mathcal{X}_{2\eta_1}$  & $\eta_1, \eta_3, \eta_4$ & $\eta_4,2\eta_3$ & $\eta_2,2\eta_4$ & $\eta_3, 2\eta_2$ & $\bullet$ & $\eta_3, 2\eta_4$ & $\eta_4, 2\eta_2$ & $\eta_2, 2\eta_3$ \\
\hline
 $\mathcal{X}_{2\eta_2}$ & $\eta_4, 2\eta_3$ & $\eta_1, \eta_3, \eta_4$ & $\eta_1, 2\eta_4$ & $\eta_1, 2\eta_3$ & $\eta_3, 2\eta_4$ & $\bullet$ & $\eta_4, 2\eta_1$ & $\eta_3, 2\eta_1$ \\
\hline
$\mathcal{X}_{2\eta_3}$ & $\eta_2, 2\eta_4$ & $\eta_1, 2\eta_4$ & $\eta_1, \eta_2, \eta_4$ & $\eta_1, 2\eta_2$ & $\eta_4, 2\eta_2$ & $\eta_4, 2\eta_1$ & $\bullet$ & $\eta_2, 2\eta_1$ \\
\hline
$\mathcal{X}_{2\eta_4}$ & $\eta_3, 2\eta_2$ & $\eta_1, 2\eta_3$ & $\eta_1, 2\eta_2$ & $\eta_1, \eta_2, \eta_3$ & $\eta_2, 2\eta_3$ & $\eta_3, 2\eta_1$ & $\eta_2, 2\eta_1$ & $\bullet$ \\
\hline
\end{tabular}
\end{equation}
Table \ref{tavola int biellittiche} has to be read in the obvious way: for instance, it says that the scheme-theoretical intersection of $\mathcal{X}_{\eta_1}$ and $\mathcal{X}_{\eta_2}$ is $N_{\eta_3}+2N_{\eta_4}$, the scheme-theoretical intersection of $\mathcal{X}_{\eta_1}$ and $\mathcal{X}_{2\eta_1}$ is $N_{\eta_2}+N_{\eta_3}+N_{\eta_4}$, and so on. \\
We are now ready to prove the main result of this section.
\begin{theorem} \label{4basepoints}
The linear system $|\mathfrak{D}_0|$ contains at most the four base points 
\begin{displaymath}
0+\eta_1+2\eta_1, \ldots, 0+\eta_4+2 \eta_4.
\end{displaymath}
\end{theorem}
\begin{proof}
For any non-zero $3-$torsion point $\eta_i$, let us consider the reducible divisor $\mathcal{Y}_{\eta_i}:=\mathcal{X}_{\eta_i}+D_{2\eta_i}$. Since $\mathcal{X}_{\eta_i}\in |-K_{E(3)}+\eta_i|$, it follows $\mathcal{Y}_{\eta_i} \in |\mathfrak{D}_0|$. Then it is sufficient to show that, set-theoretically:
\begin{equation} \label{4basepointsB} 
\bigcap_{i=1}^8 \mathcal{Y}_{\eta_i}=\{0+\eta_k+2\eta_k\; | \; k=1,2,3,4 \}.
\end{equation}
Looking at Table \ref{tavola int biellittiche}, we find that the intersection of the eight divisors $\mathcal{X}_{\eta_i}$ is empty, as well as any intersection of seven of them. On the other hand, the intersection of four or more divisors algebraically equivalent to $D$ is empty too. For any subset $I\subset \{\eta_i\; | \; i \neq 0 \}$, let us define:
\begin{displaymath}
\mathcal{X}_I:=\bigcap_{\eta \in I}\mathcal{X}_{\eta}, \quad D_I:=\bigcap_{\eta \in I}D_{\eta}, \quad 2I:=\{2 \eta\; | \; \eta \in I \},
\end{displaymath}
whereas $|I|$ denotes as usual the number of elements in $I$. Therefore, according to the above remark we have:
\begin{displaymath} 
\bigcap_{i=1}^8 \mathcal{Y}_{\eta_i}= \bigcap_{i=1}^8(\mathcal{X}_{\eta_i} \cup D_{2 \eta_i})= \mathscr{A} \cup \mathscr{B}
\end{displaymath}
where
\begin{displaymath} 
\begin{split}
\mathscr{A}&:=\bigcup_{\substack{ |I|=6,\; |J|=2\\|I \cup J|=8 \\ (I \cup 2I) \cap J= \emptyset    }} \bigg( \mathcal{X}_I \cap  D_J \bigg),\\   \mathscr{B}&:=\bigcup_{\substack{ |I|=5,\; |J|=3\\|I \cup J|=8 \\ (I \cup 2I) \cap J= \emptyset    }} \bigg( \mathcal{X}_I \cap  D_J \bigg).
\end{split}
\end{displaymath}
Looking again at table \ref{tavola int biellittiche}, it follows that $\mathscr{A}$ contains exactly four non empty terms, namely:
\begin{displaymath}
\begin{split}
\mathcal{X}_{\eta_1} \cap \mathcal{X}_{2\eta_1} \cap \mathcal{X}_{\eta_2} \cap \mathcal{X}_{2\eta_2} \cap \mathcal{X}_{\eta_3} \cap \mathcal{X}_{2\eta_3} \cap D_{\eta_4} \cap D_{2\eta_4}&=N_{\eta_4} \cap D_{\eta_4} \cap D_{2\eta_4}=\{0+\eta_4+2\eta_4 \} \\
\mathcal{X}_{\eta_1} \cap \mathcal{X}_{2\eta_1} \cap \mathcal{X}_{\eta_2} \cap \mathcal{X}_{2\eta_2} \cap \mathcal{X}_{\eta_4} \cap \mathcal{X}_{2\eta_4} \cap D_{\eta_3} \cap D_{2\eta_3}&=N_{\eta_3} \cap D_{\eta_3} \cap D_{2\eta_3}=\{0+\eta_3+2 \eta_3 \} \\
\mathcal{X}_{\eta_1} \cap \mathcal{X}_{2\eta_1} \cap \mathcal{X}_{\eta_3} \cap \mathcal{X}_{2\eta_3} \cap \mathcal{X}_{\eta_4} \cap \mathcal{X}_{2\eta_4} \cap D_{\eta_2} \cap D_{2\eta_2}&=N_{\eta_2} \cap D_{\eta_2} \cap D_{2\eta_2}=\{0+\eta_2+2 \eta_2 \}\\ 
\mathcal{X}_{\eta_2} \cap \mathcal{X}_{2\eta_2} \cap \mathcal{X}_{\eta_3} \cap \mathcal{X}_{2\eta_3} \cap \mathcal{X}_{\eta_4} \cap \mathcal{X}_{2\eta_4} \cap D_{\eta_1} \cap D_{2\eta_1}&=N_{\eta_1} \cap D_{\eta_1} \cap D_{2\eta_1}=\{0+\eta_1+2 \eta_1 \}.
\end{split}
\end{displaymath}
On the other hand, the only terms in $\mathscr{B}$ that could give contribution are of the form:
\begin{displaymath}
 \mathcal{X}_{\pm \eta_i} \cap  \mathcal{X}_{\eta_j} \cap \mathcal{X}_{2\eta_j} \cap \mathcal{X}_{\eta_k} \cap \mathcal{X}_{2\eta_k} \cap D_{\pm 2\eta_i} \cap D_{\eta_l} \cap D_{2\eta_l} 
\end{displaymath}
where $(i,j,k,l)$ is a permutation of $(1,2,3,4)$. Using again table \ref{tavola int biellittiche} one sees that all the intersections of this type are empty, hence $\mathscr{B}= \emptyset$.
Therefore we proved equality (\ref{4basepointsB}), and we are done. \\
\end{proof}
\section{A remarkable element inside $|\mathfrak{D}_0|$} \label{remarkable} 
We have hitherto described some property of the linear system $|\mathfrak{D}_0|$, whose general element is a smooth surface of general type with $p_g=q=1, \; K^2=g=3$. From the point of view of the deformation theory, it would be interesting to know how kind of singular elements are contained in $|\mathfrak{D}_0|$. We already know that it contains the eight reducible divisors $\mathcal{Y}_{\eta_i}=\mathcal{X}_{\eta_i}+D_{2\eta_i}$ and the pencil $D_0+|-K_{E(3)}|$; hence it is natural to wonder if there are \emph{irreducible}, singular elements. In this section we show that the answer is \emph{yes}, by explicitly exhibiting one of them. Let us start by considering the following three curves in $E(3)$: for any $i=1,2,3$ let us define
\begin{displaymath}     
\ell_i= \{ \xi_i+x+(x \oplus \xi_i) \; |\;x \in E \}.
\end{displaymath}
\begin{proposition}
We have
\begin{displaymath}
\ell_i D=1, \quad \ell_i F=2.
\end{displaymath}
\end{proposition}
\begin{proof}
Let $u \in E$ be a general point. Set theoretically we have $\ell_i \cap D_u= \ell_i \cap D_{u \oplus \xi_i}= \{\xi_i+u+(u \oplus \xi_i) \}$; on the other hand, since $D_u$ and $D_{u \oplus \xi_i}$ intersect transversally at the point $\xi_i+u+(u \oplus \xi_i)$, we obtain $\ell_i D=1$. Let us consider now the set-theoretical intersection of $\ell_i$ with a fibre of $E(3)$, say $F_0$; this is equal to
\begin{displaymath}
\ell_i \cap F_0= \{ 0+\xi_i+\xi_i, \; \xi_1+\xi_2+\xi_3 \}.
\end{displaymath}
On the other hand, since $\ell_i$ is a smooth elliptic curve which dominates $E$, the restriction to $\ell_i$ of the projection $\beta \colon E(3) \longrightarrow E$ must be {\'e}tale; this in turn implies that $\ell_i$ intersects transversally the fibres of $\beta$, hence $\ell_iF=2$.
\end{proof}
Let us define the following surface in $E(3)$:
\begin{displaymath}
\mathcal{Y}:=\{x+(x \oplus y)+y \; | \; x,y \in E \}.
\end{displaymath}
\begin{proposition} \label{ipsilon}
The following holds:
\begin{enumerate}
\item $\mathcal{Y}$ is an irreducible element of $|\mathfrak{D}_0|$;
\item the curves $\ell_1, \; \ell_2, \; \ell_3$ are double curves of $\mathcal{Y}$ and $\mathcal{Y}$ contains no other singularities;
\item the desingularization of $\mathcal{Y}$ is ruled.
\end{enumerate}
\end{proposition}
\begin{proof}
($1$) Let $F_u$ be a general fibre of $E(3)$; if $v$ is one of the four points of $E$ such that $2v=u$, then the intersection of $\mathcal{Y}$ and $F_u$ is given by the four curves 
\begin{equation} \label{fourlines}
\begin{split}
&r_0:=\{v + x + (v \div x )\; | \; x \in E \}, \quad \textrm{and } \\
&r_i:=\{(v \oplus \xi_i) + x + (v \oplus \xi_i \div x )\; | \; x \in E \} \quad i=1,2,3.  
\end{split}
\end{equation}
Via the identification of $F_u$ with $\mathbb{P}^2$, these curves are four lines, hence $\mathcal{Y}$ is algebraically equivalent to $4D+aF$. On the other hand, an explicit computation shows that, for general $u,v \in E$, the surface $\mathcal{Y}$ intersects transversally the elliptic curve $D_u \cap D_v=\{u+v+x \; | \; x \in E \}$ at  the three points $u+v+(u\oplus v), \; u+v+(v \div u), \; u+v+(u \div v)$. Therefore $(4D+aF)D^2=3$, that is $\mathcal{Y}$ is algebraically equivalent to $4D-F$. Let us consider now the scheme-theoretical intersection $\mathcal{Y} \cap D_u$; it consists of the two curves $\{x+(x\oplus u)\; | \; x \in E\}$ and $\{x+ (u \div x) \; | \; x \in E \}$. If one identifies $D_u$ with $E(2)$ in the natural way, these curves are linearly equivalent respectively to $4h_0-2f_0$ and $f_u$, hence the intersection $\mathcal{Y} \cap D_u$ is reduced and it is linearly equivalent to $4h_0-f_u$. This implies $\mathcal{Y} \in |\mathfrak{D}_0|$. Finally, $\mathcal{Y}$ contains no divisors algebraiucally equivalent to $D$, hence it must be irreducible. \\  
($2$) Since any curve of type $D_u \cdot \mathcal{Y}$ is reducible, $\mathcal{Y}$ turns out to be singular, because the general paracanonical curve of a surface with $p_g=q=1, \; K^2=g=3$ is irreducible (see Theorem \ref{catanese-ciliberto}). In order to understand the singularities of $\mathcal{Y}$, let's look more carefully at the intersection of $\mathcal{Y}$ with $F_u$. We already remarked that such a intersection consists of four ``lines'' $r_0, \ldots, r_3$ whose explicit equations are given by (\ref{fourlines}). The points in which these lines intersect are
\begin{displaymath}
\begin{aligned}  
P_1 & =v +(v \oplus \xi_1) + \xi_1;  & P_2 & =v +(v \oplus \xi_2) + \xi_2; & P_3 & =v +(v \oplus \xi_3) + \xi_3; \\
Q_1 & =(v \oplus \xi_2)+(v \oplus \xi_3) + \xi_1; &  Q_2 & =(v \oplus \xi_1)+(v \oplus \xi_3) + \xi_2;  & Q_3 & =(v \oplus \xi_1)+(v \oplus \xi_2) + \xi_3. 
\end{aligned} 
\end{displaymath}
If $u \neq 0$ these are six distinct points, and so the four lines $r_0, \ldots, r_4$ form a complete quadrilateral in $F_u$. If $u=0$ we have only four distinct points, namely $\xi_1+\xi_2+\xi_3, \; 0+\xi_1+\xi_1, \; 0+\xi_2+\xi_2, \; 0+\xi_3+\xi_3$; in this case the quadrilateral degenerates, because there are three lines passing through the first point. On the other hand, notice that for any $i=1,2,3$ the curve $\ell_i$ intersects the fibre $F_u$ transversally at the two points $P_i, \; Q_i$ and not elsewhere; from this is follows that $\ell_1,\; \ell_2, \;\ell_3$ form the whole singular locus of $\mathcal{Y}$. They are double curves of $\mathcal{Y}$ since the singularities of $\mathcal{Y} \cap F_u$ are ordinary nodes for $u \neq 0$; moreover, $\ell_1,\; \ell_2, \;\ell_3$ meet at the point $\xi_1+\xi_2+\xi_3$ and not elsewhere. \\  
($3$) Since $\mathcal{Y}$ is covered by rational curves, its desingularization is a ruled surface.
\end{proof}
\begin{remark}
It is immediate to check that $\mathcal{Y}$ contains the four points $0+\eta_1+2\eta_1, \ldots ,0+\eta_4+2\eta_4$, cf. Theorem \ref{4basepoints}.
\end{remark}
\section{Smooth elements of $|\mathfrak{D}_0|$} \label{smoothD0}
Theorem \ref{centrale} tells us that a smooth element $X \in |\mathfrak{D}_0|$ is a  minimal surface of general type with $p_g=q=1, \; K^2=g=3$ and ample canonical class; conversely, any minimal surface with these invariants and ample canonical class is embedded in $E(3)$ as a divisor $ X \in |\mathfrak{D}_0|$ by means of a translate of the paracanonical map. In the sequel we always identify $S$ with $X$ whenever $K$ is ample. This  section is devoted to establish some numerical results which hold in this situation.
\begin{lemma} \label{notreperduepunti}
If $K$ is ample, then do not exist two points $x,y \in S$ such
that the paracanonical curves passing through $x$
coincide with the paracanonical curves passing through
$y$.
\end{lemma}
\begin{proof}
If $K$ is ample then the paracanonical map $\omega \colon
S \longrightarrow E(3)$ is an immersion, that is the paracanonical
system separates the points on $S$.
\end{proof}
\begin{lemma}\label{razionale}
If $K$ is ample, then a paracanonical divisor and an Albanese
fibre have at most one common irreducible component $A$. Moreover
in this case $A$ is a reduced, smooth rational curve with $A^2=-3, \; KA=1$.
\end{lemma}
\begin{proof}
Since $K$ is ample, the surface $S$ is a smooth divisor in $E(3)$
algebraically equivalent to $4D-F$, and the paracanonical system
$\{K\}$ is cut out on $S$ by the divisors of the system $\{D\}$.
On the other hand, we can identify a fibre of the Abel-Jacobi map 
$\beta \colon E(3) \longrightarrow E$ with
$\mathbb{P}^2$, and with this identification the divisors of
$\{D\}$ cut out on this fibre an elliptic $1-$dimensional system of
lines; then $A$ is one of these lines, so it is a reduced, smooth rational
curve; moreover $KA=1$ because two distinct lines in
$\mathbb{P}^2$ intersect in a single point.
\end{proof}
\begin{corollary} \label{unasolarazionale}
If $K$ is ample, then a paracanonical curve $C$ contains at most one rational component $A$, appearing with multiplicity one. Moreover $A$ is smooth and $A^2=-3$.
It follows that every irreducible component $B$ of $C$ with $p_a(B)=1$ is
smooth.
\end{corollary}
\begin{proof}
This is an immediate consequence of Lemma \ref{razionale}, since a
rational component of $C$ must be contained in a fibre of the
Albanese map.
\end{proof}
\begin{proposition} \label{lemmadispezzamento}
Let us suppose that $K$ is ample, and let $C$ be a reducible
paracanonical divisor on $S$.
Then it is reduced, and there are only the following possibilities for $C$:\\ \\
$(1a)$ \quad $C=A+B$ with $A^2=-2, \; KA=2, \; B^2=-3, \; KB=1, \; AB=4$;\\
$(1b)$ \quad $C=A+B$ with $A^2=0, \; KA=2, \; B^2=-1, \; KB=1, \; AB=2$;\\
$(2a)$ \quad $C=A+B_1+B_2$ with $B_i^2=A^2=-1, \; KB_i=KA=1, \\ \phantom \quad \quad \quad \quad \quad  B_1B_2=1, \; B_iA=1$;\\
$(2b)$ \quad $C=A+B_1+B_2$ with $B_i^2=-1, \; KB_i=1, \; A^2=-3, \; KA=1,\\ \phantom \quad \quad \quad \quad \quad  B_1B_2=0, \; B_iA=2$.\\ \\
Moreover, all the irreducible components of $C$ are smooth and
 intersect transversally except possibly in the case $(1b)$.
\end{proposition}
\begin{proof}
First of all, we prove that the paracanonical divisors are
reduced. In fact, let $C$ be a reducible paracanonical divisor;
since $K$ is ample and $K^2=3$, $C$ contains at most three
irreducible components. So, if $C$ is not reduced, we would have either
$C=3A$ or $C=2A+B$. The former case is impossible because we would
have $A^2=1/3$. So we have to consider only the latter case. Clearly
$KA=KB=1$, and by the Index Theorem $K^2A^2\leq(KA)^2, \; K^2B^2 \leq (KB)^2$. Then
$A^2\leq0, \; B^2\leq0$, so the genus formula implies that $A$ and
$B$ are either rational or elliptic. By Corollary \ref{unasolarazionale},
$C$ contains at most one rational component that appears with multiplicity one; then we can suppose that $A$ is not
rational, that is $A$ is elliptic with $A^2=-1$. The arithmetic
genus of $C$ is:
\begin{displaymath}
4=p_a(C)=p_a(2A)+p_a(B)+2AB-1,
\end{displaymath}
and this implies $5=p_a(B)+2AB$, that is  $p_a(B)=1, \; AB=2$.
Now we have $C=2A+B=A+(A+B)$ and $A(A+B)=1$; this is a contradiction, because $C$ is 2-connected ([Bo73, Lemma 1 p.181]), hence $C$ is reduced.\\
Now, assume that $C$ contains only two distinct connected
components, $C=A+B$.
Without loss of generality, we can suppose $KA=2, \; KB=1$; then the Index Theorem implies either $B^2=-3$ or $B^2=-1$. On the other hand we have $KB=(A+B)B$, that is  $AB=1-B^2$. Therefore $B^2=-3$ implies $AB=4$, hence $A^2=-2$, whereas $B^2=-1$ implies $AB=2$, that is $A^2=0$. These are the cases $(1a)$ and $(1b)$.\\ 
Finally  assume that $C$ contains three distinct connected
components, that is $C=A_1+A_2+A_3$. Clearly $KA_i=1$, so by Index
Theorem $(A_i)^2=-1$ or $(A_i)^2=-3$.
Then we have the following four cases:\\
\begin{displaymath}
\begin{array}{ll}
(1)\quad (A_1)^2=-1, \quad (A_2)^2=-1, \quad (A_3)^2=-1; \\
(2)\quad (A_1)^2=-1, \quad (A_2)^2=-1, \quad (A_3)^2=-3;\\
(3)\quad (A_1)^2=-1, \quad (A_2)^2=-3, \quad (A_3)^2=-3; \\
(4)\quad (A_1)^2=-3, \quad (A_2)^2=-3, \quad (A_3)^2=-3. \\
\end{array}
\end{displaymath}
Cases ($3$) and ($4$) are excluded by Corollary \ref{unasolarazionale}; if ($1$) occurs, then from the relations $(A_1+A_2+A_3)A_i=1$ we obtain $A_1A_2+A_1A_3=A_1A_2+A_2A_3=A_1A_3+A_2A_3=2$, that is $A_iA_j=1$ for $i\ne j$, and this is the case $(2a)$. In a similar way  we can see that ($2$) corresponds to case $(2b)$.\\
Since all  the irreducible components in the cases $(1a), \; (2a), \; (2b)$ are rational or elliptic, Corollary \ref{unasolarazionale} implies that they are smooth curves.\\
To complete the proof , it remains to show that all the
irreducible components of $C$ intersect transversally except in the case
$(1b)$. This is obvious in case $(2a)$; in case $(1a)$, since $B$
is a component of an Albanese fibre, a multiple intersection would
correspond to a ramification point of the morphism
$\alpha|_{A} \colon A\longrightarrow E$ and this is a contradiction because a
morphism between two elliptic curves is {\'e}tale; the case $(2b)$ is
similar.
\end{proof}
\begin{remark} \label{altropuntodivista}
It is possible to give a nice geometric interpretation of Proposition \ref{lemmadispezzamento}. Indeed, if $K$ is ample then the paracanonical map $\omega$ gives an isomorphism between $S$ and its canonical model $X \in |\mathfrak{D}_0|$. Moreover  the paracanonical system $\{K\}$ is cut out on $X$ by the surfaces of $\{D \}$. But it is also possible to "reverse" the point of view, that is we can see the paracanonical curve $C:=C_u$ as the divisor cut out by $X$ on $D_u$.  Consider  the standard identification of $D_u$ with $E(2)$, and let $f$ be the class of a fibre and $h$ be the tautological class. By Theorem \ref{centrale}, we have $C \equiv 4h-f$ as a divisor in $E(2)$; therefore, if $C$ is reducible, Theorem \ref{coomologia in E(n)} implies that we have only the following possibilities for its irreducible components:\\
$(1a)$ $C=A+B$, where $A\equiv 4h-2f, \; B\equiv f$;\\
$(1b)$ $C=A+B$, where $A\equiv 3h-f, \; B\equiv h$;\\
$(2a)$ $C=A+B_1+B_2$, where $ A\equiv 2h-f,\;B_1,B_2 \equiv h$;\\
$(2b)$ $C=A+B_1+B_2$, where $ A \equiv f, \;B_1, B_2 \equiv 2h-f$. \\ 
It is clear that these four possibilities are exactly the cases which appear in Proposition \ref{lemmadispezzamento}. We remark that the curve $A$ in case $(1a)$ belongs to $|-2K_{E(2)}|$, whereas the curve $A$ in case $(2a)$ and the curves $B_1, \; B_2$ in case $(2b)$ are algebraically but not linearly equivalent to $-K_{E(2)}$ $($see Section \ref{symmetric products}$\;)$.
\end{remark}
Now, we recall that  a fibre of $E(3)$ cuts out on $D_u \cong
E(2)$ a fibre $f$ and on $S$ an Albanese fibre $F$; an immediate
consequence of this fact and of the previous remark is the
following:
\begin{corollary} \label{intersezionealbanese}
With the same notations of Proposition \ref{lemmadispezzamento}, we have:\\
$(1a)$ $FA=4, \; FB=0$;\\
$(1b)$ $FA=3, \; FB=1$;\\
$(2a)$ $FA=2, \; FB_1=FB_2=1$;\\
$(2b)$ $FA=0, \; FB_1=FB_2=2$.
\end{corollary}
\begin{proposition} \label{comp a comune}
If $K$ is ample, then two distinct paracanonical curves $C_u$ and $C_v$
have at most one common component, which is a smooth elliptic
curve $A$ with $KA=1$.
\end{proposition}
\begin{proof}
The components of the scheme-theoretical intersection of $C_u$ and $C_v$ are reduced by
Proposition \ref{lemmadispezzamento} and they must be contained in
the intersection $D_u \cap D_v$, which is the smooth elliptic curve $A=\{ x+u+v \; | \; x \in E \}$.
\end{proof}
\begin{proposition} \label{ellittica orizz}
If $K$ is ample, then any elliptic curve $A$ which is contained in a fibre of the Albanese pencil of $S$ satisfies $KA=3$.
\end{proposition}
\begin{proof}
If one identifies the fibres of $E(3)$ with $\mathbb{P}^2$, $A$ must be a plane cubic curve, hence $KA=3$.
\end{proof}  
 
\part{Surfaces of general type with $p_g=q=1, \; K^2=3$ and a rational pencil of genus $2$.} \label{part:genus2}

\section{Reducible paracanonical divisors} \label{reducible parac}
Throughout this section, $S$ is a minimal surface of general
type with $p_g=q=1$, $K^2=g=3$ containing a pencil $|G|$ of curves of genus $2$. $|G|$ is necessarily a rational pencil, because $\{F\}$ is the only irrational pencil on $S$ (this follows from the universal property of the Albanese map). 
Let us recall the following:
\begin{proposition} \label{curve del pencil}
There exists exactly one irreducible family of minimal surfaces of
general type  with $p_g=q=1$, $K^2=3$ and a rational pencil $|G|$
of curves of genus $2$, and this family is parametrized by the
coarse moduli space of elliptic curves. Moreover:
\begin{enumerate}
\item $g=3$; 
\item the pencil $|G|$ is base point free; 
\item $|G|$ is the only genus $2$ pencil on $S$;
\item $|G|$ contains $13$ singular fibres; six of these are elliptic curves
with an ordinary double point, the other seven consist of two
smooth elliptic curves intersecting transversally at a single
point.
\end{enumerate}
\end{proposition}
\begin{proof}
The first assertion is [Xi85b, Theorem 6.5]; (1) follows from
[Xi85b, Theorem 6.5], (2) is [Ho77, Theorem 5], (3) follows from [Xi85b, Theorem 6.5] and (4) from [Xi85b, Lemma 3.11 and Theorem 3.16]. 
\end{proof}
Since the dimension of the moduli space $\mathscr{M}$ of surfaces with $p_g=q=1, \; K^2=g=3$ is $5$, it follows that the general surface with these invariants does not contain any genus $2$ pencil.
\begin{notation}
We denote by $\mathscr{M}_2$ the moduli 
space of surfaces with $p_g=q=1, \; K^2=g=3$ which contain a genus $2$ pencil.
\end{notation}
Proposition \ref{curve del pencil} tells us that $\mathscr{M}_2$ is a subvariety of $\mathscr{M}$ isomorphic to the $j-$line. The construction of $\mathscr{M}_2$ given by Xiao in [Xi85b] involves Jacobian fibrations and modular forms. The aim of Section \ref{th1} will be to explicitly describe $\mathscr{M}_2$ as a subscheme of $\mathscr{M}$ (Theorem \ref{theorem 1}); in order to do this, we have to prove some general facts about the paracanonical system of surfaces belonging to $\mathscr{M}_2$. 
Our first step is to show, by means of a numerical argument, that if $S$ contains a genus $2$ pencil, then it doesn't contain $(-2)-$curves. This will allow us to apply to $S$ the results of Section \ref{smoothD0}.
\begin{proposition} \label{ampiezzacanonico}
If $S$ contains a pencil $|G|$ of curves of genus $2$, then $K$ is
ample.
\end{proposition}
\begin{proof}
Assume by contradiction that $K$ is not ample; then $S$ contains a
$(-2)-$curve $\Gamma$. By the Index Theorem we have $(\Gamma+G)^2
K^2 \leq [(\Gamma+G)K]^2$, that is $3(-2+2\Gamma G) \leq 4$; then
$\Gamma G \leq 1$. Moreover, part $(4)$ of Proposition \ref{curve del pencil}
shows that $\Gamma$ is not a component of a reducible fibre of
$|G|$; then $\Gamma G \geq 1$, that is $\Gamma G=1$. Let now
$W_1+Z_1, \ldots ,W_7+Z_7$ be the seven reducible fibres of the
pencil $|G|$; since $\Gamma G=1$, without loss of generality we
can suppose $W_i\Gamma =0, \; Z_i\Gamma =1$ for any $i$. We
claim that the classes of the curves $W_1,\ldots ,W_7, \; \Gamma, \; K, \;G$
are independent in Num($S$). Indeed, assume a numerical
equivalence relation:
\begin{equation} \label{eq:equivalnum}
\sum_{k=1}^7 a_kW_k+b\Gamma+cK+dG \equiv0
\end{equation}
 holds, where $a_k,b,c,d$ are integers.
Multiplying (\ref{eq:equivalnum}) by $W_i, \; \Gamma, \;  K, G$ we obtain
the following relations among the coefficients:
\begin{displaymath}
\left\{ \begin{array}{ll}
-a_i+c=0 \quad (i=1, \ldots ,7) \\
-2b+d=0 \\
\sum_{k=1}^7 a_k+3c+2d=0 \\
b+2c=0 .
\end{array} \right.
\end{displaymath}
Equivalently:
\begin{displaymath}
\left\{ \begin{array}{ll}
-2b+d=0 \\
5c+d=0 \\
b+2c=0 ,
\end{array} \right.
\end{displaymath}
and since this homogeneous linear system has only the zero
solution, our claim is proved. Then: $$\textrm{rank Num}(S) \geq
10.$$ On the other hand Noether's formula gives $c_2(S)=9$, then
$h^{1,1}(S)=c_2(S)-2+4q-2p_g=9$, so: $$\textrm{rank Num}(S) \leq
h^{1,1}(S)=9,$$ which is a contradiction.
\end{proof}
\begin{lemma} \label{pencillemma1}
Let $u \in E$, $u \neq 0$ and let $G$ be a general element of $|G|$. Then $h^0(G, \; C_u|_G)=1$.
\end{lemma}
\begin{proof}
Let us consider the short exact sequence:
\begin{displaymath}
0 \longrightarrow \mathcal{O}_S (C_u-G) \longrightarrow \mathcal{O}_S(C_u)
\longrightarrow \mathcal{O}_G(C_u|_G) \rightarrow 0.
\end{displaymath}
From this, since $u\neq 0$, we obtain the following exact sequence
of cohomology groups:
\begin{displaymath}
\label{eq:genere2} 0 \longrightarrow H^0(S, \; C_u-G) \longrightarrow H^0(S, \; C_u)
\longrightarrow H^0(G,\; C_u|_G) \longrightarrow H^1(S, \; C_u-G) \longrightarrow 0.
\end{displaymath}
Now $H^2(S, \; C_u-G)=H^0(S, \; K-C_u+G)=0$. In fact let us suppose that $G' \in |K-C_u+G|$ is effective; then since $GG'=0$, $G'$ must be union of components of fibres of $|G|$.
On the other hand we have $(G')^2=0$, hence Zariski's lemma (see [BPV, p.90]) implies $G' \cong rG$ for some $r \in \mathbb{Q}$. But we have $G' \equiv G$ as well, then we obtain $G' \cong G$ and this is impossible since $u \neq 0$. \\  
Then Riemann Roch theorem gives:
\begin{displaymath}
h^0(S, \; C_u-G)-h^1(S, \; C_u-G)= \frac{1}{2}(C_u-G)(C_u-G-K)+1=0,
\end{displaymath}
and from the exact sequence (\ref{eq:genere2}) we obtain
$h^0(G, \;C_u|_G)= h^0(S, \; C_u)=1$ as we desired.
\end{proof}
\begin{corollary} \label{pencilcoroll1}
Let $G$ be a smooth curve in the genus $2$ pencil, and let $u \in
E, \; u \neq 0$. If $C_u$ cuts out on $G$ the divisor $p+q$, then
$C_{\div u}$ cuts out the divisor $\overline{p}+\overline{q}$,
where $\overline{p}, \; \overline{q}$ are the points which are
conjugates to $p, \; q$ with respect to the hyperelliptic involution of $G$.
\end{corollary}
\begin{proof}
Since $C_u+C_{\div u} \in |2K|$, $C_u+C_{\div u}$ cuts out on $G$
an element of the bicanonical series, which is composed with the
canonical $g^1_2$. Now Lemma \ref{pencillemma1} implies $p+q
\notin g^1_2$, so the result follows.
\end{proof}
\begin{corollary} \label{pencilcoroll1'}
Let $G$ be a smooth curve in the genus $2$ pencil and let $p \in
G$. If $\{ C_u,C_v,C_w \}$ are the paracanonical curves through
$p$, then $\{C_{\div u}, C_{\div v}, C_{\div w}\}$ are the
paracanonical curves through $\overline{p}$. In particular, if
$p\in C_t$, then $\overline{p} \in C_{\div t}$.
\end{corollary}
\begin{proof}
Immediate by  Corollary \ref{pencilcoroll1}.
\end{proof}
\begin{corollary} \label{pencilcoroll1a}
Let $G$ be a smooth curve in the genus $2$ pencil and $p \in G$ be
a Weierstrass point. If $p \in C_u$, then also $p \in C_{\div u}$.
\end{corollary}
\begin{proof}
Obviously we can suppose $u\neq 0$. Let $p+q$ be the divisor cut
out on $G$ by $C_u$; since $p$ is a Weierstrass point, Lemma 
\ref{pencillemma1} implies $p \neq q$. Therefore Corollary 
\ref{pencilcoroll1} implies that $C_{\div u}$
cuts out on $G$ the divisor $\overline{p}+\overline{q}=
p+\overline{q}$.
\end{proof}
\begin{corollary} \label{pencilcoroll2}
Let $G$ a smooth curve in the genus $2$ pencil. Then any
paracanonical divisor $C_{\xi_i}, \; i \neq 0$, cuts out on $G$ two
distinct Weierstrass points.
\end{corollary}
\begin{proof}
Let $p+q$ be the divisor cut out on $G$ by $C_{\xi_i}$; from
Lemma \ref{pencillemma1} we have $p+q \notin |K_G|$. Moreover,
since $2C_{\xi_i} \in |2K|$, we have $2p+2q \in |2K_G|$; but
$|2K_G|$ is composed with $|K_G|$, so we get $2p, \; 2q \in |K_G|$, and this means that $p,\; q$ are Weierstrass points of $G$. 
Finally, Lemma \ref{pencillemma1} implies $p \neq q$.
\end{proof}
The next one is a crucial step:
\begin{proposition} \label{paracanoniciriducibili}
Let us suppose that $S$ contains a genus $2$ pencil $|G|$. Then the three
paracanonical curves $C_{\xi_i}, \; i\neq 0$, are reducible.
\end{proposition}
\begin{proof}
Let us suppose that $C_{\xi_i}$ is irreducible and let $p \in C_{\xi_i}$ be a general point. Let $G_p$ be the curve of $G$ which contains $p$; since $p$ is general, $G_p$ is a smooth genus $2$ curve. Then by Corollary \ref{pencilcoroll2} we know that $p$ is a Weierstrass point of $G_p$, and since $p \in C_{\xi_i}$ Corollary \ref{pencilcoroll1a} shows that the three paracanonical curves passing through $p$ are $\{ C_{\xi_i},\;  C_u, \; C_{\div u} \}$ for some $u \in E$. Then $p \in F_{\xi_i \oplus u \div u}=F_{\xi_i}$, and since $p \in C_{\xi_i}$ is a general point this means that $C_{\xi_i}$ is a component of $F_{\xi_i}$. Using Zariski's lemma we get $(C_{\xi_i})^2 \leq 0$, which is a contradiction.
\end{proof} 
\begin{remark}
On the other hand, it will turn out that the canonical curve $K$ of $S$ is smooth $($and hyperelliptic$)$: see Proposition \ref{covering}.
\end{remark}
Proposition \ref{ampiezzacanonico} shows that if $S$ contains a genus $2$ pencil $|G|$ then $K$ is ample. Since we have just proven that the paracanonical divisors $C_{\xi_i}, \; i=1,2,3$ are reducible, each of them must fall into one of the cases described in Proposition \ref{lemmadispezzamento}. The situation is explained by the following:
\begin{proposition} \label{decomposizione}
The paracanonical divisors $C_{\xi_i}, \; i \neq 0$ are reducible of type $(2b)$ $($see Proposition \ref{lemmadispezzamento}$\;)$. More precisely we can write:
\begin{equation} \label{A,B}  
C_{\xi_i}=A_{\xi_i}+B_{ij}+B_{ik},
\end{equation}
where $FA_{\xi_i}=0, \; FB_{ij}=FB_{ik}=2$ and $(i,j,k)$ is a permutation of $(1,2,3)$. Moreover, if we identify $S$ with its canonical model $X \in |\mathfrak{D}_0|$, we have:
\begin{equation} \label{A,C}
\begin{split}
A_{\xi_i} & = \{ \xi_i+x+(\div x) \; | \; x \in E\}; \\
B_{ij} & = \{\xi_i+x+(x \oplus \xi_j) \; | \; x \in E \}; \\
B_{ik} & = \{\xi_i+x+(x \oplus \xi_k) \; | \; x \in E \}.
\end{split}
\end{equation}
\end{proposition}
\begin{proof}
Let us identify $S$ with $X \in |\mathfrak{D}_0|$.
The argument in the proof of Proposition \ref{paracanoniciriducibili} actually shows that the curves $C_{\xi_i}$ and $F_{\xi_i}$ have a common component. Let us denote it by $A_{\xi_i}$; then $A_{\xi_i}=\{\xi_i +x+ (\div x) \thinspace | \thinspace x \in E\}$.
Now, let us consider the curve $C_{\xi_1}$; it is cut out on $X$ by the
divisor $D_{\xi_1} \subset E(3)$. If we identify $D_{\xi_1}$ with
$E(2)$ in the natural way, the curve $A_{\xi_1}$ corresponds to
the fibre $f_0$, and $X$ cuts out on $E(2)$ a divisor
linearly equivalent to $4h_0-f_{\xi_1}$ (see Lemma
\ref{restringere}). Therefore if $B:=C_{\xi_1}-A_{\xi_1}$ and we define $b$ as the curve cut out by $B$ on $D_{\xi_1}$, then the following linear equivalence relation holds:
\begin{displaymath}
4h_0-f_{\xi_1} \cong f_0+b,
\end{displaymath}
that is $b \cong -2K_{E(2)}+\xi_1$. Therefore the first
example after Proposition \ref{plurianticanonici} implies that up to the identification of $D_{\xi_1}$ with $E(2)$ we have $b=T_{\xi_2}+T_{\xi_3}$, hence: \begin{displaymath}
B=\{\xi_1+x+(x\oplus \xi_2) \thinspace | \thinspace x \in E \}+\{\xi_1+x+(x\oplus \xi_3) \thinspace | \thinspace x \in E \},
\end{displaymath}
and this gives the desired decomposition of $C_{\xi_1}$. Of course
we can repeat exactly the same argument for $C_{\xi_2}$ and
$C_{\xi_3}$, and we are done.
\end{proof}
\begin{remark} \label{prerid}
Each $B_{ij}$ is a smooth elliptic curve, and we have:
\begin{displaymath}
B_{ij} \cdot B_{kl}= \left\{ \begin{array}{ll}
1 & \emph{if } \;$i=l, \; j=k$ \\
-1 & \emph{if } \; $i=k, \; j=l$ \\ 
0 & \emph{otherwise}. 
\end{array} \right.
\end{displaymath}
\end{remark}
\begin{proposition} \label{ennei}
If $X \in |\mathfrak{D}_0|$ contains a genus two pencil $|G|$, then the four elliptic curves:
\begin{displaymath}
N_i:=N_{\eta_i}=\{x+(x \oplus \eta_i)+(x \oplus 2\eta_i)\; |\; x \in E \}\quad (i\neq 0)
\end{displaymath}
are contained in $X$.
\end{proposition}
\begin{proof}
Let $i \in \{1,2,3,4\}$ be fixed. Proposition \ref{enneeta} tells us $N_iD=1, \; N_iF=3$; therefore in $E(3)$ we have $N_iX=N_i(4D-F)=1$. On the other hand, the three curves $A_{\xi_i}$ are contained in $X$ by Proposition \ref{decomposizione}, hence the three points $\xi_1+(\xi_1\oplus \eta_i) + (\xi_1 \oplus 2\eta_i), \; \xi_2+(\xi_2\oplus \eta_i) + (\xi_2 \oplus 2\eta_i), \; \xi_3+(\xi_3\oplus \eta_i) + (\xi_3 \oplus 2\eta_i)$ belong to $X$. Since these points belong to $N_i$ as well, it follows by B{\'e}zout that $N_i$ is contained in $X$.  
\end{proof}
\begin{corollary} \label{emmei}
Let $X$ as above and $i\neq 0$. Then the two paracanonical divisors $C_{\eta_i}$ and $C_{2\eta_i}$ are both reducible and they contain a common component, namely the elliptic curve
\begin{displaymath}
M_i:=\{x+\eta_i+2\eta_i \;| \; x \in E \}.
\end{displaymath}
\end{corollary}
\begin{proof}
Let $i \in \{1,2,3,4\}$ be fixed.
Since the curves $A_{\xi_1}, \; A_{\xi_2}, \; A_{\xi_3}$ are contained in $X$, the three points $\xi_1+\eta_i+2\eta_i, \; \xi_2+\eta_i+2\eta_i, \; \xi_3+\eta_i+2\eta_i$ belong to $X$. Moreover, Proposition \ref{ennei} implies that the point $0+\eta_i+2\eta_i$ belongs to $X$ as well (cf. Theorem \ref{4basepoints}). Since $C_{\eta_i} \cdot C_{2\eta_i}=3$, it follows by B{\'e}zout that $C_{\eta_i}$ and $C_{2 \eta_i}$ have a common component, which must be the curve $M_i$ according to Proposition \ref{comp a comune}.    
\end{proof}
We are now able to completely describe the reducible fibres of $|G|$ (cf. Proposition \ref{curve del pencil}).
\begin{proposition} \label{riducibili}
Let $X \in |\gD|$ be a divisor containing a genus $2$ pencil $|G|$. Then
the seven reducible fibres of $|G|$ are the following:
\begin{displaymath}
\begin{aligned}
G_1 &:= B_{23}+B_{32}; &  G_2&:= B_{13}+B_{31}; &  G_3&:= B_{12}+B_{21};\\
G_1^*&:=M_1+N_1; & G_2^*&:=M_2+N_2; & &   \\
G_3^*&:=M_3+N_3; & G_4^*&:=M_4+N_4. 
\end{aligned}
\end{displaymath}
The seven singular points of these fibres are respectively:
\begin{displaymath}
\begin{aligned}
P_1&=0+\xi_2+\xi_3; & P_2&=0+\xi_1+\xi_3; & P_3&=0+\xi_1+\xi_2;\\
P_1^*&=0+\eta_1+2\eta_1; & P_2^*&=0+\eta_2+2\eta_2; & &\\ 
P_3^*&=0+\eta_3+2\eta_3; & P_4^*&=0+\eta_4+2\eta_4. 
\end{aligned}
\end{displaymath}
\end{proposition}
\begin{proof}
By Remark \ref{prerid} we have $(B_{ij}+B_{ji})^2=0$; since moreover $G(B_{ij}+B_{ji})=0$, it follows by Zariski's lemma that $G_1,\; G_2, \; G_3$ are reducible fibres of $|G|$. Let us consider now $G_1^*:=M_1+N_1$. We have $M_1N_1=1$, and moreover $M_1^2=N_1^2=-1$, hence $(M_1+N_1)^2=0$. Since $M_1B_{23}=N_1B_{23}=M_1B_{32}=N_1B_{32}=0$, it follows $G_1(M_1+N_1)=0$ and one concludes again by using Zariski's lemma. 
\end{proof}
As a consequence, we obtain the following result.
\begin{corollary} \label{unaS}
There is exactly one smooth surface in $|\mathfrak{D}_0|$ containing a genus $2$ pencil. 
\end{corollary}
\begin{proof}
Since surfaces with $p_g=q=1$, $K^2=g=3$ and a genus $2$ pencil do exist, the linear system $|\mathfrak{D}_0|$ contains at least one such surface $X$. Suppose that there is another one, say $X'$. Then $X'$ cuts out on $X$ a divisor linearly equivalent to $4K-F_0$ which contains, according to Proposition \ref{riducibili}, the seven curves $G_1, \ldots G_3,\; G_1^*, \ldots,G_4^*$. This would imply that $4K-F_0-7G$ is an effective class of divisors in $X$, a contradiction because $(4K-F-7G)K=-6 <0$.  
\end{proof}
\begin{notation}
In the sequel of this paper we denote by $\mathcal{X}$ the unique smooth surface in $|\mathfrak{D}_0|$ containing a genus $2$ pencil.
\end{notation} 
\begin{proposition} \label{FG}
The N{\'e}ron-Severi group of $\mathcal{X}$ has rank $9$ and it is generated by the following classes:
\begin{displaymath}
K,\;G,\;B_{12},\;B_{13},\;B_{23},\;M_1,\;M_2,\;M_3, \; M_4.
\end{displaymath}
Moreover, the following linear relations hold in $\emph{Pic}(\mathcal{X})$:
\begin{enumerate}
\item $3K \cong 3G+ \sum_{i=1}^3 A_{\xi_i}$\\
\item $F_0+3G\cong 2K+ \sum_{i=1}^4M_i$ \\
\item $F_0+\sum_{i=1}^4N_i \cong 2K+G$.
\end{enumerate}
As a consequence, we have $FG=4$.
\end{proposition}
\begin{proof}
Since rk $NS(\mathcal{X})\leq h^{1,1}(\mathcal{X})=9$, it is sufficient to show that the above classes are numerically independent, and this can be done by using exactly the same argument as in proof of Proposition \ref{ampiezzacanonico}.\\
Relation ($1$) is obtained by summing up the three equations $C_{\xi_i}=A_{\xi_i}+B_{ij}+B_{ik}$, and by using Proposition \ref{riducibili}. Consider now the surface $Y:=\{x+(\div x)+y \;|\;x,y \in E \}$. One can easily see (for instance, by intersecting it with $D_0$) that $Y \cong D_0+F_0$; hence $Y$ cuts out on $\mathcal{X}$ a curve linearly equivalent to $K+F_0$. On the other hand, such a curve contains $A_{\xi_1}, \; A_{\xi_2}, \; A_{\xi_3}, \; M_1, \; M_2, \; M_3,\; M_4$; since $K(K+F_0)=K(\sum_{i=1}^3 A_{\xi_i}+\sum_{j=1}^4 M_j)$ and $K$ is ample, it follows $K+F_0 \cong \sum_{i=1}^3 A_{\xi_i}+\sum_{j=1}^4 M_j$; thus using relation ($1$) we obtain ($2$). ($3$) is a consequence of ($2$) since Proposition \ref{riducibili} implies $\sum_{i=1}^4M_i \cong 4G-\sum_{i=1}^4N_i$. Finally, multiplying the two sides of ($1$) by $F$ we get $FG=4$ and we are done.
\end{proof}

\section{The description of the locus $\mathscr{M}_2$} \label{th1}
The following is the main result of Part \ref{part:genus2}. 
\begin{theorem} \label{theorem 1}
If a divisor in $|\mathfrak{D}_0|$ contains one of the three curves $A_{\xi_i} \; (i \neq 0)$, then it contains \emph{all} the curves $A_{\xi_1}, \; A_{\xi_2}, \; A_{\xi_3}$. Moreover there exists exactly one such a divisor, namely $\mathcal{X}$. 
\end{theorem}
Theorem \ref{theorem 1} tells us that for any $i=1, \; 2, \; 3$ we have:
\begin{displaymath}
h^0(E(3), \; \mathfrak{D}_0 \otimes \mathscr{I}_{A_{\xi_i}})=h^0(E(3), \; \mathfrak{D}_0 \otimes \mathscr{I}_{A_{\xi_1}} \otimes \mathscr{I}_{A_{\xi_2}}\otimes \mathscr{I}_{A_{\xi_3}})=1,
\end{displaymath}
and the only effective divisor in the linear system $|\mathfrak{D}_0 \otimes \mathscr{I}_{A_{\xi_i}}|$ is the smooth surface $\mathcal{X}$. It is clear that, by varying the isomorphism class of $E$, we obtain a description of $\mathscr{M}_2$ as a subvariety of codimension $4$ in $\mathscr{M}$. 
\begin{proof}
Let $\xi_i$ a $2-$torsion point on $E$ different from $0$. We firstly show that we have $h^0(E(3), \; \mathfrak{D}_0 \otimes \mathscr{I}_{A_{\xi_i}})=1$. Let us make the natural identification between $D_{\xi_i}$ and
$E(2)$ (see Section \ref{symmetric products}); therefore Lemma
\ref{restringere} implies:
\begin{align}
\mathcal{O}_{E(3)}(F_0) \otimes \mathcal{O}_{D_{\xi_i}} & = \mathcal{O}_{E(2)}(f_{\xi_i}), \label{restrizioni a D 1}\\
\mathcal{O}_{E(3)}(D_0) \otimes \mathcal{O}_{D_{\xi_i}} & = \mathcal{O}_{E(2)}(h_0), \label{restrizioni a D 2}\\
\intertext{and moreover:} \label{restrizioni a D 3}
\mathcal{O}_{D_{\xi_i}}(A_{\xi_i})& = \mathcal{O}_{E(2)}(f_0).
\end{align}
Consider the exact sequence of sheaves in $E(3)$:
\begin{equation} \label{seq 3}
0 \longrightarrow \mathcal{O}(\mathfrak{D}_0-D_{\xi_i}) \longrightarrow
\mathcal{O}(\mathfrak{D}_0) \longrightarrow
\mathcal{O}(\mathfrak{D}_0)|_{D_{\xi_i}} \longrightarrow 0.
\end{equation}
Since $\xi_i$ is a $2-$torsion point different from $0$, we
obtain:
\begin{equation} \label{vanishing1}
\begin{split}
 h^0( E(3), \; \mathfrak{D}_0-D_{\xi_i}) &=h^0(E(3), \; 4D_0-F_0-D_{\xi_i}) \\ &=h^0(E(3),\; -K_{E(3)}+\xi_i)=0\end{split}
\end{equation}
(see Proposition \ref{plurianticanonici}), hence also $h^1(E(3), \;
\mathfrak{D}_0-D_{\xi_i})=0$ by Proposition
\ref{coomologia in E(n)}. It follows that there is an isomorphism:
\begin{equation} \label{coomologia1}
 H^0(E(3), \; \mathfrak{D}_0) \xrightarrow{\; \; \cong \;} H^0(D_{\xi_i}, \; \mathfrak{D}_0|_{D_{\xi_i}}).
\end{equation}
Now we have a commutative diagram:
\begin{displaymath} 
\xymatrix{
0 \longrightarrow H^0(E(3), \; \mathfrak{D}_0 \otimes \mathscr{I}_{A_{\xi_i}}) \ar[r] 
\ar[d] & H^0(E(3), \; \mathfrak{D}_0) \ar[d] \\
0 \rightarrow H^0(D_{\xi_i}, \; (\mathfrak{D}_0 \otimes \mathscr{I}_{A_{\xi_i}})|_{D_{\xi_i}}) 
\ar[r] 
& H^0(D_{\xi_i}, \; \mathfrak{D}_0|_{D_{\xi_i}}) }
\end{displaymath}
and this shows that the isomorphism (\ref{coomologia1}) induces an isomorphism:
\begin{displaymath} 
 H^0(E(3), \; \mathfrak{D}_0 \otimes \mathscr{I}_{A_{\xi_i}}) \xrightarrow{ \; \;  \cong \;} H^0(D_{\xi_i}, \; (\mathfrak{D}_0 \otimes \mathscr{I}_{A_{\xi_i}})|_{D_{\xi_i}}).
\end{displaymath}
Applying relations (\ref{restrizioni a D 1}),
(\ref{restrizioni a D 2}) and (\ref{restrizioni a D 3}) we get:
\begin{displaymath} 
\begin{split}
H^0(D_{\xi_i}, \; (\mathfrak{D}_0 \otimes
\mathscr{I}_{A_{\xi_i}})|_{D_{\xi_i}}) &= H^0(E(2), \; 4h_0-f_{\xi_i}-f_0)\\
&=H^0(E(2), \; -2K_{E(2)}+\xi_i)=\mathbb{C},
\end{split}
\end{displaymath}
that is $h^0(E(3), \; \mathfrak{D}_0 \otimes \mathscr{I}_{A_{\xi_i}})=1$ for $i=1,2,3$. By  Proposition \ref{decomposizione}, $\mathcal{X}$ contains the three curves $A_{\xi_1}, \; A_{\xi_2}, \; A_{\xi_3}$. Thus, if $X_i$ denotes the unique effective divisor inside the linear system $|\mathfrak{D}_0 \otimes \mathscr{I}_{A_{\xi_i}}|$, we obtain $X_1=X_2=X_3=\mathcal{X}$. This completes the proof.
\end{proof}
\begin{remark}
On the other hand we have $h^0(E(3), \; \gD \otimes \mathscr{I}_{A_{\xi_0}})=3$, see Proposition \ref{uguale3}.
\end{remark}
\section{Surfaces with bicanonical map of degree $2$.} \label{proof of theorem 2}
Let $S$ be a minimal surface of general type and $\sigma \colon S \longrightarrow S$ be a birational involution. Then $\sigma$ is biregular, and its fixed locus is given by a (possibly reducible) smooth curve $R'$ and isolated fixed points $p_1, \ldots, p_t$. Let $\pi \colon \hS \longrightarrow S$ be the blow-up of $S$ at $p_1, \ldots, p_t$; therefore $\sigma$ extends to an involution $\hat{\sigma}\colon \hS \longrightarrow \hS$ whose fixed locus is:
\begin{displaymath}
\hR = \hR' + \sum_{i=1}^t E_i,
\end{displaymath}
where $E_i$ is the exceptional divisor over $p_i$ and $\hR'$ is a smooth curve isomorphic to $R'$. Let $T:=S / \sigma, \; \hW:=\hS / \hat{\sigma}$ and let 
$ \psi \colon S \longrightarrow T, \;  \hat{\psi} \colon \hS \longrightarrow \hW$ be the projections onto the quotients. The surface $T$ has $t$ nodes, whereas $\hW$ is smooth and we have a commutative diagram:  
\begin{displaymath} 
\begin{CD}
\hS  @>\pi>> S\\
@VV{\hat{\psi}}V  @VV{\psi}V\\
\hW @> \rho>> T,\\
\end{CD}
\end{displaymath}
where $\rho$ is the blow-up of $T$ at the nodes. $\hat{\psi}$ is a double cover and its branch locus $\hB$ is given by:
\begin{displaymath}
\hB= \hB'+ \sum_{i=1}^t \Omega_i,
\end{displaymath}
where the $\Omega_i$'s  are $(-2)-$curves. Let $\hmL$ be the element in Pic$(\hW)$ such that $2 \hmL=\hB$ and which determines $\hat{\psi}$. Then we have $\hat{\psi}_* \mathcal{O}_{\widehat{S}}=\mathcal{O}_{\widehat{W}} \oplus \hmL^{-1}$, where $\mathcal{O}_{\widehat{S}}$ is the invariant part and $\hmL^{-1}$ is the anti-invariant part of $\hat{\psi}_*\mathcal{O}_{\widehat{S}}$ under the action of $\hat{\sigma}$. Since $\hat{\psi}$ is a double cover, the invariants of $\widehat{S}$ and $\widehat{W}$ are related in the following way (see [BPV84, p.183]):
\begin{displaymath} 
\begin{array}{ll}
 K_{\widehat{S}}^2 & = \quad 2(K_{\widehat{W}}+\hmL)^2,\\
 \chi(\mathcal{O}_{\widehat{S}}) & = \quad 2 \chi(\mathcal{O}_{\widehat{W}}) + \frac{1}{2} \hmL  \cdot (K_{\widehat{W}}+\hmL),\\
 p_g(\widehat{S}) & = \quad p_g(\widehat{W})+h^0(\widehat{W}, \; K_{\widehat{W}}+\hmL),\\
 q(\widehat{S}) & = \quad q(\widehat{W})+h^1(\widehat{W}, \; K_{\widehat{W}}+\hmL).
\end{array}
\end{displaymath}
Let us suppose now that the bicanonical map of $S$ has degree $2$. Then $S$ possess the so called \emph{bicanonical involution} $\sigma \colon S \longrightarrow S$, which exchanges the two sheets of the (generically) double cover $\phi \colon S \longrightarrow \Sigma$. Notice that by definition $\phi$ is composed with $\sigma$, that is we have a commutative diagram:
\begin{displaymath} 
\begin{CD} 
\xymatrix{
\hS \ar[rrrr]^{\hat{\phi}} \ar[rrd]_{\hat{\psi}} & & & & \Sigma \\
& & \widehat{W} \ar[rru]_{\mu} & &
} 
\end{CD}
\end{displaymath}
where the map $\mu$ is birational. 
The existence of such a diagram implies that the following equalities hold (see [CM02, Proposition 6.1]):
\begin{equation} \label{relazionirivestimento}
\begin{array}{ll}
(i) & (2K_{\widehat{W}}+\hB')^2=2K_S^2;\\
(ii) & \chi(\mathcal{O}_{\widehat{W}}(2K_{\widehat{W}}+\hmL))=0;\\
(iii) & K_{\widehat{W}} \cdot (K_{\widehat{W}}+\hmL)= \chi(\mathcal{O}_{\widehat{W}})-\chi(\mathcal{O}_S).
\end{array}
\end{equation}
If $S$ is a surface of general type with a pencil of curves of genus $2$, then the bicanonical map $\phi$ of $S$ is not birational, because $|2K|$ cuts out the general curve of the pencil a subseries of the bicanonical series, composed with the hyperelliptic involution. This is the so called \emph{standard case} for the non-birationality of the bicanonical map.
Let $d:= \deg \phi$; there is the following:
\begin{theorem}[Xiao] \label{grado=2}
Let $S$ be a minimal surface of general type with a genus $2$
pencil $|G|$. Suppose moreover $P_2(S) \geq 4$. Then $d=2$ unless
we are in one of the following cases:
\begin{enumerate}
\item $p_g=q=4, \; K^2=8$, $S$ is the product of two curves of genus
$2$; in this case $d=4$; 
\item $p_g=2, \; q=1$ or $2$, $K^2=4$, $S$
has two pencils of curves of genus $2$; in this case $d=4$; 
\item $p_g=1, \; q=0,\; K^2=4$, $S$ has the genus $2$ pencil $|G|$ and a
genus $3$ pencil $|M|$ and $GM=4$; in this case $d=4$; 
\item $p_g=1, \; q=0, \; K^2=2$ or $3$; in this case $d=2$ or $4$; 
\item some fibration with $p_g=q=2, \; K^2=8$ and $p_g=q=1, \; K^2=4$; in these
cases $d=4$.
\end{enumerate}
\end{theorem}
\begin{proof}
See [Xi85b, Theorem 5.6]
\end{proof}
\begin{remark}
The assumption $P_2(S) \geq 4$ is necessary. Indeed if we consider a surface $S$ with $p_g=q=1, \, K^2=g=2$ $($\;these surfaces are studied in $[$\emph{Ca81}$]$\;$)$, then we have $P_2(S)=3$, and the bicanonical map is a covering of degree $8$ of $\mathbb{P}^2$.
\end{remark}
\begin{corollary} \label{degbig}
The bicanonical map $\phi$ of a surface of general type with $p_g=q=1,\; K^2=g=3$ and a genus $2$ pencil $|G|$ has degree $2$.
\end{corollary}
On the other hand, we say that a surface of general type with non birational bicanonical map presents the \emph{non-standard case} if it does not contain any pencil of curves of genus $2$. There are only finitely many families of surfaces presenting the non standard case, but a complete classification is still missing. In the paper [Xi90] Xiao gave a long list of possibilities; later on several authors investigated about their real occurrence.  
In the case $p_g=q=1, \; K^2=g=3$ we can state the results of Xiao in the following Theorem \ref{Xiao list}, that we call the ``Xiao's list'':
\begin{theorem}[Xiao's list] \label{Xiao list} 
Let $S$ be a minimal surface of general type with $p_g=q=1,\; K^2=g=3$ and bicanonical map of degree $2$, presenting the non-standard case. Then only the following possibilities can occur:
\begin{itemize} 
\item[$(A)$] $\hW$ is a ruled surface. In this case $\hW$ must be \emph{rational} and $S$ contains a pencil of hyperelliptic curve of genus $3$, $5$ or $7$.
\item[$(B1)$] $\hW$ is an Enriques surface;
\item[$(B2)$] $\hW$ is a regular surface with geometric genus $1$ and Kodaira dimension $1$.
\end{itemize}  
\end{theorem}
The aim of this section is to improve this result: indeed we will show that cases $(A)$ and $(B1)$ actually don't occur (Proposition \ref{due casi grado 2}). We don't know whether case $(B2)$ occurs or not.  \\
Let $j \colon E(3) \longrightarrow E(3)$ be the automorphism given by
$j(x+y+z)= (\div x)+(\div y) + (\div z)$.
\begin{definition}
We say that a subvariety $X \subset E(3)$ is \emph{symmetric} if
$j(X)=X$.
\end{definition}
In the remainder of this section $S$ always denotes a minimal surface of general type with $p_g=q=1, \; K^2=g=3$ and bicanonical map of degree $2$, and $X \in |\mathfrak{D}_0|$ is its canonical model. Notice that we don't make the assumption that the canonical class $K$ of $S$ is ample, so $S$ \emph{a priori} can contain $(-2)-$curves, whose contraction gives rise to some rational double points on $X$. The bicanonical involution $\sigma \colon S \longrightarrow S$ of $S$ induces a biregular involution $\sigma^c \colon X \longrightarrow X$  of $X$. The following lemma shows that the behavior of $\sigma^c$ is very simple:
\begin{lemma} \label{involuzione}
$X$ is a symmetric divisor in $E(3)$ and the involution $\sigma^c: X \longrightarrow X$ coincides with the restriction of $j$ to $X$.
\end{lemma}
\begin{proof}
It is sufficient to show that, if $\{C_u, \; C_v, \; C_w \}$ are the paracanonical curves passing through a general point $x$ of $S$, then $\{ C_{\div u}, \; C_{\div v}, \; C_{\div w} \}$ are the paracanonical curves passing through the point $y= \sigma(x)$. Since any automorphism of $S$ sends paracanonical curves into paracanonical curves, there exists an automorphism $\sigma$ of $E$ such that $\sigma(C_u)=C_{\sigma(u)}$. This implies that the paracanonical curves passing through $y$ are $\{ C_{\sigma(u)}, \; C_{\sigma(v)}, \; C_{\sigma(w)} \}$. On the other hand $\sigma$ is the bicanonical involution of $S$ and $C_t + C_{\div t} \in |2K|$ for any $t \in E$, then we must have either  $\sigma(t)=t$ or $\sigma(t)= \div t$ for any $t \in E$. The first case is impossible, otherwise $\sigma$ would be the identity map; it follows $\sigma(t)= \div t$ and we are done.
\end{proof} 
Lemma \ref{involuzione} implies the existence of a commutative diagram: 
\begin{displaymath}
\begin{CD}
S  @>\omega>> E(3)\\
@VV{\sigma}V  @VV{j}V\\
S @>\omega>> E(3).\\
\end{CD}
\end{displaymath}
Notice that the fixed locus of the involution $j$ is given by four
smooth curves $A_{\xi_0}, \ldots , A_{\xi_3}$ and four isolated
points $\rho_0, \ldots, \rho_3$, where:
\begin{displaymath} 
\begin{split}
A_{\xi_i}&:=\{\xi_i+x+(\div x) \thinspace | \thinspace x \in E \}; \\
\rho_i &:=\xi_j+\xi_k+\xi_l,
\end{split}
\end{displaymath}
and $(i,j,k,l)$ is a permutation of $(0,1,2,3)$. 
Observe moreover that both $A_{\xi_i}$ and $\rho_i$ belong to $F_{\xi_i}$ for every $i$. It follows as soon the following:
\begin{corollary} \label{luogo fisso}
The fixed curve of the involution $\sigma^c \colon X \longrightarrow X$ is
contained in the disjoint union: 
\begin{displaymath}
A_{\xi_0} \sqcup \ldots \sqcup A_{\xi_3}.
\end{displaymath} 
\end{corollary}
We remark explicitly that the $A_{\xi_i}$'s are four smooth rational curves in $E(3)$ such that $DA_{\xi_i}=1$; hence if $A_{\xi_i}$ is contained in $X$, then its strict transform in $S$ is a smooth, rational curve with self-intersection $-3$; with a little abuse of notation, we denote it again by $A_{\xi_i}$. By Corollary \ref{luogo fisso} the fixed curve $R'$ of the bicanonical involution $\sigma \colon S \longrightarrow S$ is composed by some of the $A_{\xi_i}$'s plus possibly some $(-2)-$curves.\\
Now we can rule out the cases $(A)$ and $(B1)$ of Xiao's list.
\begin{proposition} \label{due casi grado 2}
Let $S$ be a minimal surface of general type with $p_g=q=1, \; K^2=g=3$. Suppose that the bicanonical map $\phi$ has degree $2$. Then either $\widehat{W}$ is ruled, and in this case $S$ contains a genus $2$ pencil, or $\widehat{W}$ is an elliptic surface of Kodaira dimension $1$ with $p_g(\widehat{W})=1, \; q(\widehat{W})=0$. If the latter possibility occurs, then $S$ presents the non-standard case, and in particular $S$ belongs to case $($B2$\;)$ of Xiao's list.
\end{proposition} 
\begin{proof}
Let $n, \; 0 \leq n \leq 4$, be the number of components of
type $A_{\xi_i}$ in  the curve $R'$, $h$ be the number of $(-2)-$curves in $R'$ and $t$ the number of isolated
fixed points of $\sigma$. Therefore the curve $B'$ is composed of $n$
disjoint smooth rational curves $\Theta_1, \ldots, \Theta_n$ such
that $\Theta_i^2=-6, \; K_{\widehat{W}}\Theta_i=4$ and of $h$ disjoint smooth
rational curves $\Gamma_1, \ldots, \Gamma_h$ such that
$\Gamma_i^2=-4, \; K_{\widehat{W}}\Gamma_i=2$. Using the notation introduced at the beginning of this section, we can write:
\begin{displaymath} 
2 \mathcal{\widehat{L}}=\hB=\sum_{i=1}^n \Theta_i+\sum_{i=1}^h \Gamma_i+ \sum_{i=1}^t
\Omega_i,
\end{displaymath}
and an easy calculation shows that the following equalities hold:
\begin{equation} \label{relazionidiramazione}
\begin{array}{ll}
\mathcal{\widehat{L}}^2 & = -\frac{1}{2}(3n+2h+t);\\
K_{\widehat{W}}\widehat{\mathcal{L}} & = 2n+h;\\
\mathcal{\widehat{L}} \hB' & =-3n-2h;\\
(\hB')^2 & =-6n-4h;\\
K_{\widehat{W}} \hB' & =4n+2h.
\end{array}
\end{equation}
Substituting (\ref{relazionidiramazione}) in
(\ref{relazionirivestimento}) we obtain equations relating $K_{\widehat{W}}^2, \; \chi(\mathcal{O}_{\widehat{W}}), \;  n, \; h$ and $t$:
\begin{equation} \label{relazionirivestimento2}
\begin{array}{ll}
(i') & 2K_{\widehat{W}}^2+5n+2h-3=0;\\
(ii') & 4K_{\widehat{W}}^2+9n+4h-t+4 \chi(\mathcal{O}_{\widehat{W}})=0;\\
(iii') & n-t +8 \chi(\mathcal{O}_{\hW})-4=0.\\
\end{array}
\end{equation}
Eliminating $K_{\widehat{W}}$ from $(i')$ and $(ii')$ in
(\ref{relazionirivestimento2}) we can write:
\begin{equation} \label{relazionirivestimento3}
\begin{array}{ll}
(i'') & -t-n+ 4\chi(\mathcal{O}_{\widehat{W}})+6=0;\\
(ii'') & -t+n +8\chi(\mathcal{O}_{\widehat{W}})-4=0.\\
\end{array}
\end{equation}
Finally, we subtract $(ii'')$ from $(i'')$ in (\ref{relazionirivestimento3}) obtaining a
single relation between $n$ and $\chi(\mathcal{O}_{\widehat{W}})$:
\begin{equation} \label{relazionirivestimento4}
\chi(\mathcal{O}_{\widehat{W}}) =\frac{1}{2}(5-n).
\end{equation}
From (\ref{relazionirivestimento4}), it follows that
there are only the two following  numerical possibilities for
$n,t$ and the invariants of the surface $\hW$:
\begin{equation} \label{possibilitanumeriche}
\begin{array}{cccc}
(a)\quad n=3, & t=7, & \chi(\mathcal{O}_{\widehat{W}})=1, & K_{\widehat{W}}^2=-(6+h);\\
(b)\quad n=1, & t=13, & \chi(\mathcal{O}_{\widehat{W}})=2, & K_{\widehat{W}}^2=-(1+h).\\
\end{array}
\end{equation}
Suppose that case ($a$) occurs. Since $n=3$, the canonical model $X$ of $S$ contains at least one curve $A_{\xi_i}$ with $i \neq 0$; therefore by Theorem \ref{theorem 1} we obtain $S=X= \mathcal{X}$, that is $S$ contains a genus $2$ pencil; moreover  Proposition \ref{ampiezzacanonico} implies $h=0$, hence $K_{\widehat{W}}^2=-6$. \\
Let us consider now case ($b)$. Since $n=1$, Theorem \ref{theorem 1} implies that the canonical model $X$ of $S$ contains the curve $A_{\xi_0}$ and that $S$ belongs to the non-standard case. We get $\chi(\mathcal{O}_{\widehat{W}})=2$, that is $\widehat{W}$ is neither rational nor Enriques. It follows that $S$ is of type $(B2)$ of Xiao's list.  
\end{proof}
\begin{remark} \label{seven fixed}
In case $($a$\;)$ the $7$ isolated fixed points of $\sigma$ coincide with the $7$ singular points of the reducible fibres of the genus $2$ pencil $|G|$.
\end{remark} 
If a surface of type ($b$) exists, then lemma \ref{involuzione} and the above discussion show that its canonical model is a symmetric element inside $|\mathfrak{D}_0 \otimes \mathscr{I}_{ A_{\xi_0}}|$. Contrary to $|\mathfrak{D}_0 \otimes \mathscr{I}_{A_{\xi_i}}|, \; i \neq 0$, which contains the only divisor $\mathcal{X}$, the linear system $|\mathfrak{D}_0 \otimes A_{\xi_0}|$ contains many elements. Indeed we can prove the following
\begin{proposition} \label{uguale3}
We have
\begin{displaymath}
h^0(E(3),\; \mathfrak{D}_0 \otimes \mathscr{I}_{A_{\xi_0}})=3.
\end{displaymath}
Moreover, the general element of $|\mathfrak{D}_0 \otimes \mathscr{I}_{A_{\xi_0}}|$ is smooth.
\end{proposition}
\begin{proof}
Let us consider the surface $\mathcal{Y}$ described in Section \ref{remarkable}. Since $\mathcal{Y}$ contains $A_{\xi_0}$, it follows by Proposition \ref{ipsilon} that it belongs to the linear system  $|\mathfrak{D}_0 \otimes \mathscr{I}_{A_{\xi_0}}|$. It follows that  $|\mathfrak{D}_0 \otimes \mathscr{I}_{A_{\xi_0}}|$ contains the linear system spanned by the pencil $D_0+|-K_{E(3)}|$, having $D_0$ as fixed part, and by the irreducible surface $\mathcal{Y}$; this  implies $h^0(E(3),\; \mathfrak{D}_0 \otimes \mathscr{I}_{A_{\xi_0}}) \geq 3$. On the other hand, it is immediate to check that the surface $\mathcal{X}$ does not contain the curve $A_{\xi_0}$, but it contains some points of it, namely the four points $P_1^*, \ldots, P_4^*$ (see Proposition \ref{riducibili}); this shows that the codimension of $H^0(E(3),\; \mathfrak{D}_0 \otimes \mathscr{I}_{A_{\xi_0}})$ in $H^0(E(3),\; \mathfrak{D}_0)$ is at least $2$, that is $h^0(E(3),\; \mathfrak{D}_0 \otimes \mathscr{I}_{A_{\xi_0}}) \leq 3$. Hence  $h^0(E(3),\; \mathfrak{D}_0 \otimes \mathscr{I}_{A_{\xi_0}})= 3$. \\
The points of $E(3)$ which are singular for all the members of the pencil $D_0+|-K_{E(3)}|$ are given by the intersection of the curve $\Gamma=\{x +(x \oplus \xi_1)+ (x\oplus \xi_2) \; | \; x \in E \}$, base locus of $|-K_{E(3)}|$, with $D_0$; this yields the three points $0+\xi_1+\xi_2, \; 0+\xi_1+\xi_3, \; 0+\xi_2+\xi_3$. Therefore the smoothness of the general element of $|\mathfrak{D}_0 \otimes \mathscr{I}_{A_{\xi_0}}|$ will follow by Bertini's theorem if we are able to show that these points don't belong to Sing($\mathcal{Y}$). But this is clear because Proposition \ref{ipsilon} shows that Sing($\mathcal{Y}$)$=\ell_1 \cup \ell_2 \cup \ell_3$, hence Sing($\mathcal{Y}$)$\cap D_0$ consists of the four points $\xi_1+\xi_2+\xi_3, \; 0+\xi_1+\xi_1, \; 0+\xi_2+\xi_2, \; 0+\xi_3+\xi_3$. \\This completes the proof.  
\end{proof}
We don't know how many symmetric elements exist in $|\mathfrak{D}_0 \otimes \mathscr{I}_{A_{\xi_0}}|$; the only example we have is the surface $\mathcal{Y}$ but it is too singular: actually it is not a surface of general type, see Proposition \ref{ipsilon}.

\section{A double cover construction} \label{double cov}
In this section and in the next one $S$ denotes afresh a surface with $p_g=q=1, \; K^2=g=3$ containing a genus $2$ pencil $|G|$. We want to give a description of it `` {\`a} la Campedelli'', that is as a double covering of a rational surface. It will turn out that this description is actually the description of the double cover associated with the hyperelliptic involution of $S$, which coincides with the bicanonical involution (see Proposition \ref{Hirzebruch}). We will prove moreover that, although the paracanonical curves $C_{\xi_1}, \; C_{\xi_2}, \; C_{\xi_3}$ of $S$ are reducible by Proposition \ref{paracanoniciriducibili}, the canonical curve $K$ is actually a smooth, hyperelliptic curve of genus $4$ (see Proposition \ref{Hirzebruch}). In Part \ref{part:bicanonicalmap} of the paper this fact will be a crucial point in the proof of Theorem \ref{main2}. \\
Let us denote by $G_i=W_i+Z_i$, $i=1, \ldots, 7$ the seven reducible curves of $|G|$ and let $P_i=W_i \cap Z_i$. Let $G$ be a
general curve of $|G|$; [Bo73, Theorem A  p.177] says
$h^1(S, \; -G)=0$, hence we can write down the following exact sequence:
\begin{displaymath} 
0 \longrightarrow H^0(S, \; K) \longrightarrow H^0(S, \; K+G) \longrightarrow
H^0(G, \; \omega_G) \longrightarrow H^1(S, \; K) \longrightarrow 0.
\end{displaymath}
This implies $h^0(S, \; K+G)=2$, and from the sequence:
\begin{displaymath} 
0 \longrightarrow H^0(S, \; K+G) \longrightarrow H^0(S, \; K+2G) \longrightarrow H^0(G, \; \omega_G) \longrightarrow 0
\end{displaymath}
we obtain $h^0(S, \; K+2G)=4$. Since $(K+2G)|_G \cong \omega_G$, $P_1, \ldots
,P_7$ are base points of $|K+2G|$; we claim that there are no
other base points. Indeed, observe that all the base points of
$|K+2G|$ belong to the unique canonical curve $K$ of $S$ because $|2G|$ is base
point free; hence, $(K+2G)K=7$ means that there are at most seven
base points, and this implies our claim. Let  $\pi \colon \widehat{S} \longrightarrow S$ be the blow-up of $S$ at the points
$P_1,\ldots,P_7$ and let $E_i$ be the exceptional divisor over $P_i$; in the sequel, the symbol $''\thinspace \thinspace \widehat{} \thinspace  \thinspace ''$  denotes the strict transform of a curve with respect to this blow-up. The linear system $|\pi^*(K+2G)-\sum_{i=1}^7
E_i|$ has dimension $3$ and it is base-point free; since
$(\pi^*(K+2G)-\sum_{i=1}^7 E_i)^2=4$, it follows that
$|\pi^*(K+2G)-\sum_{i=1}^7 E_i|$ induces a degree $2$ morphism:
\begin{displaymath}
\hat{f}_2 \colon \widehat{S} \longrightarrow \mathbb{P}^3,
\end{displaymath}
whose image $Q$ is a nondegenerate quadric surface $Q$; 
on the other hand, if $\langle \; t \; \rangle$ is a basis for $H^0(S, \; K)$ and $\langle \; s_1, s_2 \; \rangle$ is a basis for $H^0(S, \; G)$, we can complete $\langle \; s_1^2t, \; s_2^2t, \; s_1s_2t \; \rangle$ to a basis of $H^0(K+2G)$, and this shows that $Q$ is a quadric cone. The curves of the genus $2$ pencil $|G|$ are mapped $2:1$ to the lines of the cone; moreover, since $(\pi^*(K+2G)-\sum_{i=1}^7 E_i)(\pi^*K-\sum_{i=1}^7 E_i)=0$, the strict transform $\widehat{K}$ of $K$ is contracted to the vertex $v$.\\
Let us consider the linear system $|K+3G|$; exactly as above, we
obtain that its dimension is $5$ and it has exactly $7$ base
points at $P_1,\ldots, P_7$; therefore $(\pi^*(K+3G)-\sum_{i=1}^7
E_i)^2=8$ implies that the linear system
$|\pi^*(K+3G)-\sum_{i=1}^7 E_i|$  induces a morphism $\hat{f}_3$
of degree $2$ onto a $\mathbb{F}_2$ surface embedded in
$\mathbb{P}^5$. In other words, we have the following commutative
diagram:
\begin{equation} \label{diagramma}
\begin{CD}
\widehat{S}  @>id>> \widehat{S}\\
@VV{\hat{f}_3}V  @VV{\hat{f}_2}V\\
\mathbb{F}_2  @>p>>  Q,\\
\end{CD}
\end{equation}
where $p$ is the blow-up of $Q$ at the vertex $v$; we remark that from diagram  (\ref{diagramma}) it follows that the image via $\hat{f}_3$ of the curve $\widehat{K}$ is the unique $(-2)-$ section $C_0$ of $\mathbb{F}_2$.\\
Now we want to describe the branch locus $B$ of $\hat{f}_3$; 
if $B \cong 2(aL+bC_0)$, then Hurwitz formula applied to the double
covering $\hat{f}_3|_G \colon G \longrightarrow L$ implies:
$$(aL+bC_0)L=3,$$
that is $b=3$ and therefore:
\begin{equation} \label{diramazione}
B \cong 2aL+6C_0.
\end{equation}
Furthermore, we have
$$\pi^*(G_i)=\widehat{W}_i+\widehat{Z}_i+2E_i,\quad i=1, \ldots,
7$$
where $\widehat{W}_i:=\pi^*(W_i)-E_i$ and $\widehat{Z}_i:=\pi^*(Z_i)-E_i$ are elliptic curves with self-intersection $-2$ in $\widehat{S}$.\\
Let $L_i$ be the image via $\hat{f}_3$ of the curve $E_i$; $L_i$ is a fibre of $\mathbb{F}_2$ such that $\hat{f}_3^*(L_i)=\widehat{W}_i+\widehat{Z}_i+2E_i$. This implies that there are at least seven different fibres in the branch locus of $\hat{f}_3$; thus, $B=B^{\sharp}+\sum_{i=1}^7 L_i$, where $B^{\sharp}$ is a curve 
linearly equivalent to $(2a-7)L+6C_0$.\\
It is clear that $\hat{f}_3$ cannot contract horizontal curves respect to the genus $2$ pencil $|\widehat{G}|$ in $\widehat{S}$, because the ruling of $\mathbb{F}_2$ is base-point free; then if there are contracted curves, they must be irreducible components of some curve of $|\widehat{G}|$. But therefore the relations $(\pi^*(K+3G)-\sum_{i=1}^7 E_i)\widehat{W}_k= (\pi^*(K+3G)-\sum_{i=1}^7 E_i)\widehat{Z}_k=0$ and $(\pi^*(K+3G)-\sum_{i=1}^7 E_i)E_k=1$ imply that $\f$ contracts exactly the $14$ curves $\widehat{W}_i, \; \widehat{Z}_i$. It follows that the only non negligible singularities of the branch locus $B$ are $14$ ordinary quadruple points, lying on the fibres $L_i$; then the only non negligible singularities of the curve $B^{\sharp}$ are $14$ ordinary triple points $x_i, \; y_i \in L_i$, $1 \leq i \leq 7$. On the other hand, $B$ doesn't contain non-negligible singularities, since $K$ is ample by Proposition \ref{ampiezzacanonico}.\\
Now we can calculate the integer $a$ in (\ref{diramazione});
indeed by [BPV, p.183] we have:
\begin{displaymath}
K_{\widehat{S}}^2=2(K_{\mathbb{F}_2}+(aL+3C_0))^2-28,
\end{displaymath}
and this means $a=11$, that is $B \cong 22L+6C_0$ and $B^{\sharp} \cong 15L+6C_0$.\\
Let us keep the notations of Section \ref{proof of theorem 2}, and let us consider the morphism $\hat{\psi} \colon \hS \longrightarrow \hW$. Notice that
that part ($a$) of relations (\ref{possibilitanumeriche}) implies
that $\hW$ is a smooth rational ruled surface such that $K_{\hW}^2= -6$,
that is $\hW$ is a minimal Hirzebruch surface blown-up 14 times. There are
$14$ $(-1)-$curves in $\hW$, which are images via $\hat{\psi}$ of the curves
$\widehat{W}_i, \; \widehat{Z}_i$. If we contract these $(-1)-$curves,
we obtain a morphism $\widehat{S} \longrightarrow
\mathbb{F}_e$, where $e$ is a nonnegative integer and
$\mathbb{F}_e$ is the minimal Hirzebruch surface containing a
section with self-intersection $-e$; with abuse of notation, we denote such a 
 morphism again by $\hat{\psi}$. The only curves contracted by
$\hat{\psi}$ are $\widehat{W}_i, \; \widehat{Z}_i$. The next proposition shows that the map $\hat{f}_3$ and the map $\hat{\psi}$ are essentially the same. 
\begin{proposition} \label{Hirzebruch}
We have $e=2$ and moreover there exists an automorphism $\phi$ of
$\mathbb{F}_2$ such that $\phi \circ \hat{\psi}= \hat{f}_3$.
\end{proposition}
\begin{proof}
Note that both the maps $\hat{\psi}$ and $\hat{f}_3$ send the
genus $2$ pencil of $\widehat{S}$ to the ruling of the
corresponding Hirzebruch surface. Let $x \in \mathbb{F}_e$ be a
general point and $\{y_1,y_2\}:=\hat{\psi}^{-1}(x)$. Then
$y_1+y_2$ is an element of the $g_2^1$ of a genus $2$ curve $G$ of
$\widehat{S}$; since the $g_2^1$ on $G$ is unique, it follows
$\hat{f}_3(y_1)=\hat{f}_3(y_2)=z \in \mathbb{F}_2$. The map $\phi \colon
\mathbb{F}_e \longrightarrow \mathbb{F}_2$ which associates $z$ to $x$
is clearly birational; moreover, since $\hat{\psi}$ and
$\hat{f}_3$ contract exactly the same curves, $\phi$ is actually
an isomorphism, hence $e=2$ and we obtain a commutative diagram:
\begin{displaymath}
\begin{CD}
\widehat{S}  @>id>> \widehat{S}\\
@VV{\hat{\psi}}V  @VV{\hat{f}_3}V\\
\mathbb{F}_2 @>\phi>> \mathbb{F}_2 .\\
\end{CD}
\end{displaymath}
This completes the proof.
\end{proof}
Now, let us give a more precise description of the
curve $B^{\sharp}$. Proposition \ref{Hirzebruch} shows that we can identify the map $\hat{f}_3 \colon \widehat{S} \longrightarrow
\mathbb{F}_2$ with $\hat{\psi}$; in particular, the branch locus of $\hat{f}_3$ is given by the strict transforms
$\widehat{A}_{\xi_1}, \;  \widehat{A}_{\xi_2}, \; \widehat{A}_{\xi_3}$ of
the curves $A_{\xi_1}, \; A_{\xi_2}, \; A_{\xi_3}$ and by the seven
exceptional divisors $E_1, \ldots, E_7$. The curves $E_1, \ldots,
E_7$ are mapped by $\hat{f}_3$ to the seven fibres $L_1, \ldots,
L_7$ contained in the branch locus of $\hat{f}_3$; it follows that
$B^{\sharp}$ is the image via $\hat{f}_3$ of the three curves
$A_{\xi_1}, \; A_{\xi_2}, \;  A_{\xi_3}$. Hence $B^{\sharp}$ is reducible
and we can write:
\begin{displaymath}
B^{\sharp}=B^{\sharp}_1+ B^{\sharp}_2+B^{\sharp}_3,
\end{displaymath}
where $B^{\sharp}_i=\hat{f}_3(\widehat{A}_{\xi_i})$. Note that $B^{\sharp}_i$ is birational to $\widehat{A}_{\xi_i}$, hence the geometric genus of $B^{\sharp}_i$ is $0$. Moreover $\widehat{A}_{\xi_i}G=2$ implies $B^{\sharp}_i L=2$, that is $B^{\sharp}_i \cong 2C_0+a_iL$ for some $a_i \in \mathbb{Z}$.\\
Since our $\mathbb{F}_2$ is embedded in $\mathbb{P}^5$, we have
$\mathcal{O}_{\mathbb{F}_2}(1)=C_0+3L$; if $H$ is a general
hyperplane section of $\mathbb{F}_2$, $H$ does not intersect the
points coming from the contraction of
$\widehat{W}_j, \; \widehat{Z}_j$, hence $(\hat{f}_3)^*(H)$ does not
intersect $\widehat{W}_j, \; \widehat{Z}_j$. Since
$B^{\sharp}_i$ is in the branch locus of $\hat{f}_3$, we obtain:
\begin{displaymath}
\begin{split}
(C_0+3L)(2C_0+a_iL)&= HB^{\sharp}_i= \hat{f}_3^{\; *}(H) \cdot
\hat{f}_3^{\; *}(B^{\sharp}_i)\\ &=(\pi^*(K+3G)-\sum_{j=1}^7
E_j)\widehat{A}_{\xi_i}=7,
\end{split}
\end{displaymath}
and this implies $a_i=5$, that is $B^{\sharp}_i \cong 2C_0+5L$ as
we expected.\\
We can describe the curve $B^{\sharp}$ if we are able to describe the curves $B^{\sharp}_i$. Without loss of generality, we can identify $S$ with the surface $\mathcal{X} \in |\mathfrak{D}_0|$; using the notation of Proposition \ref{riducibili}, let $E_i$ be the exceptional divisor over $P_i$, and let $E_j^*$ be the exceptional divisor over $P_j^*$; let $L_i$ be the image in $\mathbb{F}_2$ of $E_i$, and $L_j^*$ the image of $E_j^*$. Of course, $L_1, \ldots, L_3, \; L_1^*, \ldots, L_4^*$ are the seven fibres of $\mathbb{F}_2$ contained in the branch locus of the map $\hat{f}_3$. Finally, let us denote by $x_{ij}$ the point in $\mathbb{F}_2$ coming from the contraction of the curve $\widehat{B}_{ij}$, and let $m_k, \; n_k$ be the points coming from the contraction of the curves $\widehat{M}_k, \; \widehat{N}_k$ respectively. Notice that:   
\begin{displaymath}
\begin{split}
x_{23}, \; x_{32} & \in L_1;\\
x_{13}, \; x_{31} & \in L_2;\\
x_{12}, \; x_{21} & \in L_3;\\
m_k, \; n_k & \in L_k^*. 
\end{split}
\end{displaymath}
Let us consider the curve $B^{\sharp}_1$. From relations (\ref{A,B}) and (\ref{A,C}) we obtain:
\begin{equation} \label{intersect}
\begin{split}
A_{\xi_1}B_{12}=A_{\xi_1}B_{13}&=2;\\
A_{\xi_1}B_{23}=A_{\xi_1}B_{32}&=1;\\
A_{\xi_1}B_{21}=A_{\xi_1}B_{31}&=0.\\
\end{split}
\end{equation}
Moreover we know by Proposition \ref{lemmadispezzamento} that all these curves intersect transversally. Therefore (\ref{intersect}) implies that $B^{\sharp}_1$ has two ordinary double points at $x_{12}, \; x_{13}$, it pass simply through  $x_{23}, \; x_{32}$ and it does not contain $x_{21}, \; x_{31}$. Moreover, since $B^{\sharp}_1 \cong 2C_0+5L$, the genus formula gives $p_a(B^{\sharp}_1)=2$; the geometric genus of $B^{\sharp}_1$ being $0$, it follows that $x_{12}$ and $x_{13}$ are the only singular points of $B^{\sharp}_1$.\\
Finally, we have:
\begin{displaymath}
1=C_{\xi_1}M_k=(A_{\xi_1}+B_{12}+B_{13})M_k=A_{\xi_1}M_k,
\end{displaymath}
because $B_{12}$ is a component of $G_3$, $B_{13}$ is a component of $G_2$ and $M_k$ is a component of $G_k^*$; in the same way we have also $A_{\xi_1}N_k=1$, and this implies that $B_1^{\sharp}$ pass simply through  the points $m_k, \; n_k$ for any $k=1, \ldots, 4$.  
Of course, with the same argument we can give a similar description of the curves $B^{\sharp}_2$ and $B^{\sharp}_3$. \\
Now, let us consider the canonical curve $K$ of $S$. Since $B^{\sharp}C_0=1$, the intersection of $B^{\sharp}$ with $C_0$ is transversal; on the other hand, $C_0$ intersects transversally any fibre and does not contain the singular points of $B^{\sharp}$, because the curve $\widehat{K}$ does not meet the curves $\widehat{M}_k, \; \widehat{N}_k$; this implies that the intersection of $C_0$ with the branch locus $B$ is transversal, that is the canonical curve $K$ is the double covering of $C_0$ branched along the $10$ distinct points of the set $C_0 \cap B$. Therefore $K$ is \emph{smooth} and \emph{hyperelliptic}. We can summarize all these considerations in the following: 
\begin{proposition} \label{covering}
If $S$ contains a pencil $|G|$ of curves of genus $2$, then
$\widehat{S}$ is the canonical resolution of the double cover of
$\mathbb{F}_2$ branched along the curve:
\begin{displaymath}
B=B^{\sharp}+\sum_{i=1}^7 L_i,
\end{displaymath}
where the $L_i$'s are distinct fibres and $B^{\sharp}$ is a curve 
linearly equivalent to $6C_0+15L$. Besides, $B^{\sharp}$ has $2$ ordinary triple points on each fibre $L_i$ and no other singularities. More precisely, we have:
\begin{displaymath}
B^{\sharp}=B^{\sharp}_1+ B^{\sharp}_2+B^{\sharp}_3,
\end{displaymath}
where $B^{\sharp}_i$ is an irreducible curve linearly equivalent to $2C_0+5L$ containing two ordinary double points.\\
The multiplicities of the curves $B^{\sharp}_i$ at the points
$x_{ij}, \; m_k, \; n_k$ are as in the following table:
\begin{displaymath} 
\begin{tabular}{|c|c|c|c|c|c|c|c|c|}
\hline
 {} & $x_{21}$ & $x_{31}$ & $x_{12}$ & $x_{32}$ & $x_{13}$ & $x_{23}$ & $m_k, \thinspace 1\leq k\leq 4$ & $n_k, \thinspace 1\leq k\leq 4$ \\
\hline
$B^{\sharp}_1$ & $0$ & $0$ & $2$ & $1$ & $2$ & $1$ & $1$ & $1$ \\
\hline
$B^{\sharp}_2$ & $2$ & $1$ & $0$ & $0$ & $1$ & $2$ & $1$ & $1$ \\
\hline
$B^{\sharp}_3$ & $1$ & $2$ & $1$ & $2$ & $0$ & $0$ & $1$ & $1$ \\
\hline
\end{tabular}
\end{displaymath}
The canonical curve $K$ of $S$ is smooth, hyperelliptic of genus
$4$.
\end{proposition}
\begin{corollary} \label{par hyp}
If $S$ contains a genus $2$ pencil $|G|$, then the general
paracanonical curve $C_u$ of $S$ is smooth, hyperelliptic of genus $4$.
\end{corollary}
\begin{proof}
The last statement of Proposition \ref{covering} says that the canonical curve $K$ of $S$ 
is smooth; hence there exists an open neighborhood $E^0$ of $0$ in $E$ such that, for any $u \in E^0$, $C_u$ is smooth.  Moreover, $|G|$ cuts out on
$C_u$ a $g^1_2$, and this implies that $C_u$ is hyperelliptic (of
genus 4).
\end{proof}
Let $\{ \widehat{F} \}$ be the Albanese pencil and $|\widehat{G}|$ the genus $2$ pencil of the surface $\widehat{S}$, whose general element are of course isomorphic to $F$ and $G$ respectively. We are interested to the image of $\{ \widehat{F} \}$ on $\mathbb{F}_2$ via the generically finite morphism $\hat{f}_3$. Let $\widehat{F}$ be a general element of $\{ \widehat{F} \}$ and $P=\hat{f}_{3}(\widehat{F})$; let $aC_0+bL$ be the class of $P$ in $\mathbb{F}_2$.
We first remark that the map $\hat{f}_3|_{\widehat{F}} \colon \widehat{F} \longrightarrow P$ is birational. Indeed $\hat{f}_3|_{\widehat{G}}$ is the canonical map of $\widehat{G}$, which identifies conjugate points with respect to the hyperelliptic involution. It follows from Corollary \ref{pencilcoroll1'} that two conjugates points of $\widehat{G}$ lie on two opposite fibres of $\{ \widehat{F} \}$, therefore $\hat{f}_3$ identifies two opposite fibres of the Albanese map of $\widehat{S}$; this means that $\hat{f}_3$ restricted to a general Albanese fibre is birational. We know that $\widehat{G}\widehat{F}=\widehat{K}\widehat{F}=4$; therefore the integers $a,b$ satisfy:
\begin{displaymath}
L(aC_0+bL)=4, \quad C_0(aC_0+bL)=4,
\end{displaymath}
that is $P\cong 4C_0+12L$. Again, we can identify $S$ with $\mathcal{X}$, and use the notation of Proposition \ref{riducibili}. Therefore we have:
\begin{displaymath}
\begin{aligned}
\widehat{F}\widehat{B}_{ij}&=2 &i,j&=1,2,3, \quad i\neq j;\\
\widehat{F}\widehat{M}_k&=1 &  k&=1,2,3,4;\\
\widehat{F}\widehat{N}_k&=3 & k&=1,2,3,4.
\end{aligned}
\end{displaymath}
Hence we obtained the following:
\begin{proposition} \label{immaginealbanese}
The image of the Albanese pencil $|F|$ of $S$ via the map
$\hat{f}_3$ is a $1-$dimensional linear system $\mathcal{F}
\subset |4C_0+12L|$. The base points of $\mathcal{F}$ are the
points $x_{ij}, m_k, n_k$ and the multiplicities of the general
curve of $\mathcal{F}$ at the base points are given in the
following table:
\begin{displaymath} 
\begin{tabular}{|c|c|c|c|c|c|c|c|c|}
\hline
 {} & $x_{21}$ & $x_{31}$ & $x_{12}$ & $x_{32}$ & $x_{13}$ & $x_{23}$ & $m_k,\thinspace 1\leq k\leq 4$ & $n_k, \thinspace 1 \leq k \leq 4$ \\
\hline
$\mathcal{F}$ & $2$ & $2$ & $2$ & $2$ & $2$ & $2$ & $1$ & $3$ \\
\hline
\end{tabular}
\end{displaymath}
\end{proposition}
Now we wish to describe the linear system $\Lambda$ in $\mathbb{F}_2$, image of the bicanonical system $|2K|$ of $S$ via the map $\hat{f}_3$. Since $|2K|$ has no base points, the same holds for the linear system $\pi^*|2K|$ in $\widehat{S}$. Let $\Gamma$ be a general element of $\pi^*|2K|$ and let $\Delta=\hat{f}_3(\Gamma)$ be the image of $\Gamma$ in $\mathbb{F}_2$; the fact that the bicanonical map of $S$ restricted to $G$ has degree $2$ implies that $\hat{f}_3$ restricted to $\Gamma$ has degree $2$, that is $\hat{f}_3|_{\Gamma} \colon \Gamma \longrightarrow \Delta$ is a double covering.
Suppose $\Delta \cong aC_0+bL$; since:
\begin{displaymath}
\Delta L=\frac{1}{2} \hat{f}_3^*(\Delta)\hat{f}_3^*(L)
=\frac{1}{2} \Gamma \widehat{G}=2,
\end{displaymath}
we obtain $a=2$.\\
Let $G=W+Z$ be a reducible fibre of $G$; since $\pi^*|2K|$
cuts out a $g^1_2$ on $\widehat{W}$, it follows that the two
intersection points of $\Gamma$ with $\widehat{W}$ are conjugate
points with respect to the involution induced over $\Gamma$ by the
map $\hat{f}_3$, and the same holds for the two intersection
points of $\Gamma$ with $\widehat{Z}$. This implies that the
curve $\Delta$ is smooth at the $14$ singular points of
$B^{\sharp}$. Since these are of course the only base points of
the linear system $\Lambda$, Bertini's theorem says that the
general curve of $\Lambda$ is smooth. On the other hand,
$\widehat{A}_{\xi_1}$ is a component of the ramification divisor 
of  $\hat{f}_3$, so
we can write:
\begin{displaymath}
B_1^{\sharp} \Delta= \widehat{A}_{\xi_1}\Gamma+14=16.
\end{displaymath}
This means:
\begin{displaymath}
(2C_0+5L)(2C_0+bL)=16,
\end{displaymath}
that is $b=7$.\\
Hence we have obtained the following result, which is a particular case of a result of Xiao (see [Xi85b, p.14]):
\begin{proposition} \label{immagine del bic}
The image in $\mathbb{F}_2$ of the bicanonical system $|2K|$ of
$S$ is the linear system $\Lambda$ given by the curves of
$|2C_0+7L|$ having simple base points at the $14$ points $x_{ij},
m_k, n_k$. 

\end{proposition}
\begin{remark}
Since $C_u +C_{\div u} \in |2K|$ for any $u \in E$, it
follows that the two paracanonical curves $C_u$ and $C_{\div u}$
are identified via the map $\phi$ to the same element of the linear
system $\Lambda$.
\end{remark}

\section{The bicanonical image} \label{bicanonical image}
Let $S$ be a minimal surface with $p_g=q=1, \; K^2=g=3$ and a genus $2$ pencil $|G|$, and let $\phi \colon S \longrightarrow \mathbb{P}^3$ be its bicanonical map; Corollary \ref{degbig} states that $\phi$ is a generically finite morphism of degree $2$. The aim of this section is to give a description of the
bicanonical image $\Sigma$ of $S$, which is a rational surface of degree $6$ in $\mathbb{P}^3$. In the sequel we keep the notation of Section \ref{double cov}, unless explicitly stated.
Let $f \colon \widehat{S} \longrightarrow S$ be the blow-up of $S$ at the seven points $P_1, \ldots, P_3,P_1^*, \ldots, P_4^*$, and $g \colon \hF \longrightarrow \F$ the blow-up of $\F$ at the $14$ points $x_{ij},m_k,n_k$. The morphism $\f \colon \hS \longrightarrow \F$ induces a rational map $\psi \colon S \dashrightarrow \F$ which is not defined at $P_1, \ldots, P_3, P_1^*, \ldots, P_4^*$; the indeterminacy of $\psi^{-1}$ is given by the 14 points $x_{ij},  m_k, n_k$.\\
We have the following commutative diagram:
\begin{displaymath} 
\begin{CD}
\hS  @>f>> S\\
@VV{\hat{\psi}}V  @VV{\psi}V\\
\hF @>g>> \F.\\
\end{CD}
\end{displaymath}
Proposition \ref{Hirzebruch} implies  $g \circ \hat{\psi}= \f$.
It is clear that $\phi$ factors through $\psi$, 
hence we have a rational map $\mu \colon \F
\dashrightarrow \Sigma$, which is induced by the linear system
$\Lambda$ described in Proposition \ref{immagine del bic}. Let
$\E_1, \ldots, \E_{14}$ be the exceptional divisors in $\hF$;
Proposition \ref{immagine del bic} says that the complete linear
system $|\La|:=|g^* \Lambda - \sum_{i=1}^{14} \E_i|$ is base
point free, hence it induces a birational morphism:
\begin{displaymath}
\hmu= \mu \circ g \colon \hF \longrightarrow \Sigma.
\end{displaymath}
In the sequel, we denote by $\La$ a general (hence smooth) curve
of $|\La|$, and by $\hL$ a general fibre of $\hF$; moreover, we
denote by $\hL_i+\E_i+\E_{7+i}$, $i=1, \ldots, 7$ the seven fibres
of $\hF$ containing the exceptional divisors.
\begin{lemma} \label{hmu}
The birational morphism $\hmu$ maps  the
exceptional divisors $\E_i, \; \E_{7+i}$ isomorphically onto lines $\lambda_i, \; \lambda_{7+i}$ and contracts the seven curves $\hL_i$; moreover it maps $\hL$ and the
curves $\widehat{B}_i^{\sharp}$ isomorphically onto conics.
\end{lemma}
\begin{proof}
The first three statements are immediate because $ \La
\E_i=\La\E_{7+i}=1, \; \La \hL_i=0, \; \La \hL=2$; the last one is a consequence of the fact that the curves $\widehat{B}_i^{\sharp}$ are in the branch locus
of the map $\hat{\psi}$ and $\widehat{B}_i^{\sharp} \La=2$.
\end{proof}
Let $H$ be a general
hyperplane section of $\Sigma$. The geometric genus of $H$ is equal
to the genus of the smooth curve $\La$; since $\La$ is isomorphic to $\Lambda$, we obtain:
\begin{displaymath} 
g(H)=\frac{1}{2} (2C_0+7L)3L+1=4.
\end{displaymath}
On the other hand, since the arithmetic genus of a plane sextic
 is $10$, it follows that any hyperplane section of $\Sigma$
is singular, and this in turn implies that $\Sigma$ has some
multiple curve. It is no difficult to describe this multiple
curve. Indeed, a divisor on $\F$ linearly equivalent to $4C_0+10L$ is
ample  by [Ha77, Chapter V, Corollary 2.18], so
$h^1(2C_0+6L)=h^1(-(4C_0+10L))=0$ by Kodaira vanishing
theorem. Then Riemann-Roch formula gives:
\begin{displaymath}
h^0(2C_0+6L)=\frac{1}{2}(2C_0+6L)(4C_0+10L)+1=15,
\end{displaymath}
that is the complete linear system $|2C_0+6L|$ has dimension $14$.
Therefore there exists at least one effective divisor $\Theta \in
|2C_0+6L|$ containing the $14$ points $x_{ij}, m_k, n_k$; let
$\Th$ be its strict transform in $\hF$. For any fibre $\hL$ in the
ruling of $\hF$,
$\mu(\Th+\hL)$ is an hyperplane section of $\Sigma$; since the image of $\hL$ is a conic, it follows that $\mu(\Th) \subset \Sigma$ is a quadruple line $\Gamma$. Hence the general hyperplane section $H$ of $\Sigma$ has a quadruple point, and since we showed that its geometric genus is $4$, this quadruple point is the only singularity of $H$; it follows that $\Gamma$ is the only multiple curve of $\Sigma$.\\
Proposition \ref{immagine del bic} says that every fibre $L_i$ $(i=1, \ldots ,7)$ contains two base points of $\Lambda$; so $\Theta L=2$ implies that the intersection of $\Theta$ with $L_i$ at these base points is transversal. As a consequence, we have that $\widehat{\Theta}$ and $\widehat{L}_i$ are disjoint for every $i=1, \ldots, 7$. So  Lemma \ref{hmu} says that $\Sigma$ contains seven isolated ordinary double points $Q_1, \ldots, Q_7$ coming from the contraction of the $(-2)-$curves $\hL_i$ via the birational morphism $\hat{\mu}$. It is easy to see that these are the only curves contracted by $\hmu$, and since $\hmu$ is birational this implies that $Q_1, \ldots, Q_7$ are the only isolated singularities of $\Sigma$.
\begin{remark}
Since $|2K|$ cuts out on the general curve of $|G|$ a
$g^2_4$, the exact sequence:
\begin{displaymath}
0 \longrightarrow H^0(S, 2K-G) \longrightarrow H^0(S, 2K) \longrightarrow
H^0(G, 2K|_G)
\end{displaymath}
implies $h^0(S, 2K-G)=1$, that is there exists just one effective
divisor $\Omega \in |2K-G|$. It is clear that $\phi(\Omega)=
\Gamma$, and since $2K(2K-G)=8$ and $\Gamma$ is a line, the
map $\phi|_{\Omega} \colon \Omega \longrightarrow \Gamma$ has degree $8$.
\end{remark}
Finally, we show that $\lambda_1, \ldots, \lambda_{14}$ and the
quadruple curve $\Gamma$ are the only lines on $\Sigma$. In order
to do this, we need the following:
\begin{lemma} \label{KA=1}
If $A \subset S$ is a curve such that $KA=1$, then one of the following 
cases holds:
\begin{enumerate}
\item $A=A_{\xi_i}$ for some $1 \leq i \leq3$; \item $A$ is one of
the $14$ elliptic curves contained in the reducible fibres of
$|G|$.
\end{enumerate}
\end{lemma}
\begin{proof}
The ampleness of $K$ implies that $A$ is irreducible, and the Index
Theorem says $A^2=-3$ or $A^2=-1$. Hence either $p_a(A)=0$ or
$p_a(A)=1$, so $A$ is smooth by Corollary \ref{unasolarazionale}. 
On the other hand, Proposition \ref{FG} yields:
\begin{displaymath}
3=3KA=3GA+\sum_{i=1}^3 A_{\xi_i}A.
\end{displaymath}
Suppose that $A$ is different from $A_{\xi_i}$ for all $1 \leq i
\leq 3$; then we have $\sum_{i=1}^3 A_{\xi_i}A \geq 0$, and since
$G$ moves in a pencil without fixed points we have only two possibilities:
\begin{enumerate}
\item[$(*)$]  $GA=0, \quad \sum_{i=1}^3 A_{\xi_i}A=3$; 
\item[$(**)$]  $GA=1,\quad \sum_{i=1}^3 A_{\xi_i}A=0$.
\end{enumerate}
If ($*$) holds, $GA=0$ implies that $A$ is a component of a reducible fibre of $|G|$, that is we are in case ($2$) of our lemma.\\
Now we suppose that ($**$) holds, and we get a contradiction.
Since $|G|$ is a rational pencil and $GA=1$, all the points of $A$
are linearly equivalent, so $A$ is smooth rational and $A^2=-3$;
moreover $GA=1$ again says that $A$ does not meet the isolated
points of the bicanonical involution $\sigma$ (because they are the
singular points of the reducible fibres of $|G|$) and
$\sum_{i=1}^3 A_{\xi_i}A=0$ shows that $A$ does not meet the fixed
curves. This implies that the involution $\sigma$ does not have fixed points on  $A$, and since $A \cong \mathbb{P}^1$ this implies $\sigma^*(\sigma(A))=A+A^*$, where $A^* \cong A$ and $AA^*=0$. Since $\sum_{i=1}^3
A_{\xi_i}A=0$ implies $\sum_{i=1}^3 A_{\xi_i}A^*=0$, we obtain that
$A_{\xi_1},A_{\xi_2},A_{\xi_3}, A, A^*$ are five disjoint rational
curves on $S$ with self-intersection $-3$. This contradicts
inequality ($6$) of [Mi84], and we are done.
\end{proof}
Now, suppose that $\lambda$ is a line in $\Sigma$, different from $\Gamma$; the general point of $\lambda$ is a smooth point for $\Sigma$, because $\Gamma$ is the only multiple curve of $\Sigma$, then if we consider the maximal divisor $A \subset S$ such that $\phi(A)=\lambda$, we have: 
\begin{displaymath}
KA=\frac{1}{2}(2K)A=H\lambda=1, 
\end{displaymath}
so $A$ is an irreducible curve on $S$ with $KA=1$, and Lemma \ref{KA=1} implies that either $A=A_{\xi_i}$ for some $1 \leq i \leq 3$ or $A$ is one of the components of the reducible fibres of $|G|$. The former case cannot hold, because Lemma \ref{hmu} shows that the images of the curves $A_{\xi_i}$ are conics; hence the latter case holds, that is $\lambda=\lambda_i$ for some $1 \leq i \leq 14$. Therefore we have the following:
\begin{proposition} \label{imag bic}
The bicanonical image $\Sigma$ of $S$ is a surface of degree $6$ in $\mathbb{P}^3$ such that its only singularities are a quadruple line $\Gamma$ and $7$ isolated ordinary double points $Q_1, \ldots, Q_7$.\\
$\Gamma$ is the image on $\Sigma$ of the only effective divisor in the linear system $|2K-G|$, whereas the double points are the images of the $7$ isolated fixed points of the bicanonical involution $\sigma$.\\
Besides the quadruple line $\Gamma$, $\Sigma$ contains exactly $14$ lines $\lambda_1, \ldots,\lambda_{14}$, which are the images of the components of the
reducible fibres of $|G|$.
\end{proposition}

\part{The bicanonical map of surfaces with $p_g=q=1, \; K^2=g=3$} \label{part:bicanonicalmap}

\section{Statement of the result} \label{result}
Let $S$ be a minimal surface of general type with $p_g=q=1$, $d$  the degree of its bicanonical map $\phi \colon S \longrightarrow \mathbb{P}^{K^2}$ and $\Sigma= \phi(S)$ the bicanonical image. A result of Xiao (see [Xi90, Proposition 5]) states that if $p_g=q=1, \; K^2 \geq 5$ and $S$ presents the non-standard case, then $d=1,2 \; \textrm{or }4$. However, this result does not apply if one wants to investigate the bicanonical map of surfaces with $p_g=q=1, \; K^2=3$. Therefore in this case, since $(2K_S)^2=12$, we have \emph{a priori} the following possibilities:\begin{equation} \label{deg poss}
\begin{split}
d=1, \quad & \textrm{deg}\; \Sigma=12; \\
d=2, \quad & \textrm{deg}\; \Sigma=6; \\
d=3, \quad & \textrm{deg}\; \Sigma=4; \\
d=4, \quad & \textrm{deg}\; \Sigma=3; \\
d=6, \quad & \textrm{deg}\; \Sigma=2.
\end{split} 
\end{equation}
The main result of Part \ref{part:bicanonicalmap} is Theorem \ref{main3} below, which says that \emph{generically} only the first case in (\ref{deg poss}) occurs.
\begin{theorem} \label{main3}
There exists a subset $\mathscr{M}^0$ of the moduli space $\mathscr{M}$ of surfaces with $p_g=q=1, \; K^2=g=3$ such that:
\begin{itemize}
\item $\mathscr{M}^0$ is dense in $\mathscr{M}$ with respect to the Zariski topology; 
\item if $[S] \in \mathscr{M}^0$, then the bicanonical map of $S$ is a birational morphism. 
\end{itemize}
\end{theorem}
Sections \ref{statement of main2} and \ref{proof main2} deal with the proof of Theorem \ref{main3}

\section{A subset of $|\mathfrak{D}_0|$ } \label{statement of main2}

Let as usual $\gD$ be a divisor in $E(3)$ linearly equivalent to $4D_0-F_0$.
\begin{proposition} \label{genlisci}
Let $u \in E$ be  a point which is not one of the eight nonzero $3-$torsion points of $E$. For a general choice of $X$ in $|\gD|$, the paracanonical divisor $C_u$ of $X$ $($i.e. the scheme-theoretical intersection  $D_u \cdot X \subset X$$)$ is a smooth, irreducible curve.
\end{proposition}
\begin{proof}
We consider separately the two cases $u \neq 0$ and $u=0$.
\begin{itemize}
\item $u\neq 0$.\\
In this case $H^0(E(3), \; \gD-D_u)=H^0(E(3), -K_{E(3)}-u)=0$ by Proposition \ref{plurianticanonici}, then $H^1(E(3), \; \gD-D_u)=0$ by Theorem \ref{coomologia in E(n)}. Hence, by considering the cohomology of the exact sequence:
\begin{displaymath}
0 \longrightarrow \mO_{E(3)}(\gD-D_u) \longrightarrow \mO_{E(3)}(\gD) \longrightarrow \mO_{E(3)}(\gD)|_{D_u} \longrightarrow 0,
\end{displaymath}
we obtain an isomorphism: 
\begin{displaymath}
H^0(E(3), \;  \gD) \xrightarrow{ \; \; \cong \;} H^0(D_u, \; \gD |_{D_u}),
\end{displaymath}
that is $|\gD|$ cuts out over $D_u$ a $complete$ linear system $|\mathcal{L}|$. Via the standard identification of $D_u$ with $E(2)$, $\mathcal{L}$ is associated with a line bundle algebraically equivalent to $4h-f$, where $f$ is a class of a fibre and $h$ is the tautological class. Thus $|\mathcal{L}|$ is base point free by [CaCi93, Theorem 1.18], and Bertini's theorem allows us to conclude that the general element of $|\mathcal{L}|$ is smooth and irreducible. \\
\item $u=0$.\\
Look at the restriction map:
\begin{displaymath}
H^0(E(3), \; \gD) \xrightarrow{\; \; r \;} H^0(D_0, \; \gD|_{D_0}).
\end{displaymath}
The image of $r$ gives a linear system $\mathcal{L}$ on $D_0$ which is without fixed components, because $|\gD|$ has at most a finite number of base points.
Moreover the last statement of Proposition \ref{covering} says that the only surface $\mathcal{X}$ in  $|\gD|$ containing a genus $2$ pencil cuts out on $D_0$ a smooth, irreducible curve. Of course this curve is an element of $\mathcal{L}$, therefore we conclude again by using Bertini's theorem.
\end{itemize}
\end{proof}
\begin{notation}
We denote by $\mathfrak{B}$ the subset of $|\gD|$ consisting of divisors $X$ enjoing the following properties: 
\begin{itemize}
\item $X$ is smooth $($that is the corresponding $S$ has ample canonical class$)$;
\item the curves $\Co:=K_X, \; \Cone, \; \Ctwo, \; \Cthree$ are smooth.
\end{itemize}
\end{notation}
Notice that if $X \in \mathfrak{B}$ then the general paracanonical curve of $X$ is smooth (hence irreducible, because the paracanonical divisors are $2-$connected).
\begin{corollary} \label{insieme denso}
$\mathfrak{B}$ is a dense subset of $|\gD|$.
\end{corollary}
\begin{proof}
Since the general elements of $|\gD|$ is smooth, the result follows from Proposition \ref{genlisci}.
\end{proof}
Our goal is to show the following:
\begin{theorem} \label{main2}
If $X \in \mathfrak{B}$, then the bicanonical map of $X$ is a birational morphism.
\end{theorem}
Let us define $\mathscr{M}^0$ as the subset of $\mathscr{M}$ parameterizing isomorphism classes of surfaces $S \in \mathfrak{B}$. It  is clear that Theorem \ref{main3} follows as soon by using Corollary \ref{insieme denso} and Theorem \ref{main2}. Therefore we need to prove Theorem \ref{main2}, and this is done in the next section.

\section{Proof of Theorem \ref{main2}} \label{proof main2}

We prove Theorem \ref{main2} by showing that, if $X \in \mathfrak{B}$, then none of cases $d=2, \; 3, \; 4, \; 6$ occurs.
\subsection{The case $d=2$} \label{d=2}
\begin{proposition} \label{non2}
If $X \in \mathfrak{B}$, then $d=2$ does not occur.
\end{proposition}
\begin{proof}
Suppose that $d=2$ occurs. Therefore the proof of Proposition \ref{due casi grado 2} tells us that there are two possibilities:
\begin{itemize}
\item $X$ contains the three curves $A_{\xi_1}, \; A_{\xi_2}, \; A_{\xi_3}$ but not the curve $A_{\xi_0}$. In this case $X=\mathcal{X}$, i.e. $X$ contains a genus $2$ pencil;
\item $X$ contains the curve $A_{\xi_0}$, but none of the curves $A_{\xi_i}, \thinspace 1 \leq i \leq3$. If this happens, then we are in the non$-$standard case.
\end{itemize}
Recall now that $A_{\xi_i}$ is a component of $C_{\xi_i}$ for $i=0,1,2,3$. Thus in the former case the three paracanonical divisors $\Cone, \; \Ctwo, \; \Cthree$ are reducible, whereas in the latter case the canonical divisor $K_X=\Co$ is reducible. In both cases, $X$ is not contained in $\mathfrak{B}$ and we are done.
\end{proof}
\subsection{The case $d=3$} \label{d=3}
\begin{proposition} \label{finite rette}
Suppose $d=3$, and let $k \in \mathbb{Z}$ be an odd positive integer. Then $\Sigma$ contains at most finitely many irreducible plane curves of degree $k$. In particular, $\Sigma$ contains at most finitely many lines and finitely many plane cubics.
\end{proposition}
\begin{proof}
Suppose that $\Sigma$ contains infinitely many plane curves of degree $k$. Let $T$ be one of them such that $T$ is contained neither in the singular locus of $\Sigma$ nor in the branch locus of $\phi$. If $\Theta$ is the pull-back of $T$ in $S$, then we have:
\begin{displaymath}
K \Theta=\frac{1}{2}(2K) \Theta= \frac{3}{2}k,
\end{displaymath}
a contradiction because $k$ is odd.
\end{proof}
\begin{corollary} \label{curva multipla}
If $d=3$ and $L \subset \Sigma$ is a line, then $\Sigma$ contains a multiple curve $\Gamma$ such that $L \subset \Gamma$.
\end{corollary} 
\begin{proof}
Suppose that the general point of $L$ is smooth for $\Sigma$; since $\Sigma$ is not ruled in lines by Proposition \ref{finite rette}, the pencil of planes through $L$ cuts out a pencil of cubics in $\Sigma$. But this contradicts again Proposition \ref{finite rette}.  
\end{proof}
\begin{corollary} \label{infinite coniche}
If $d=3$, then $\Sigma$ contains an infinite number of conics.
\end{corollary}
\begin{proof}
Let $u \in E$ be a general point. We know by Theorem \ref{catanese-ciliberto} that the general paracanonical divisor $C_u \in S$ is irreducible and moreover that $C_u+C_{\div u} \in |2K|$. Then the image in $\Sigma$ of $C_u+C_{\div u}$ is a plane curve of degree $4$. Since $(2K)C_u=6$ and $4$ does not divide $6$, it follows that $\phi(C_u)$ and $\phi(C_{\div u})$ are distinct, irreducible curves. Therefore Proposition \ref{finite rette} implies that, for general $u$, $\phi \colon C_u \longrightarrow \phi(C_u)$ is triple cover of a smooth conic, and we are done. 
\end{proof}
Let us denote by $H(u, \div u)$ the plane of $\mathbb{P}^3$ that cuts out on $\Sigma$ the plane section $C_u + C_{\div u}$; therefore $\mathcal{H}:=\{H(u,\div u) \; | \; u \in E\}$ is an algebraic, one-dimensional system of planes in $\mathbb{P}^3$ parametrized by the rational curve $\Delta=\{C_u+C_{\div u}\}_{u \in E} \subset |2K|$. Proposition 3.3 of [CaCi91] tells us that $\Delta$ is an irreducible curve of degree $3$ in the $3-$dimensional projective space $|2K|$, in particular $\Delta$ is not a line. Thus we obtain:
\begin{proposition} \label{nopencil}
The algebraic system $\mathcal{H}$ is not a pencil, hence the planes of $\mathcal{H}$ don't contain a common straight line.
\end{proposition}
By Corollary \ref{infinite coniche} the surface $\Sigma$ is a quartic surface that contains an algebraic system of plane sections which are reducible into pairs of conics.  
The study of such surfaces was made by Kummer in [Ku]; a good account of this subject is Chapter VII of [Je16]. In order to state the result of Kummer, we need the following remark: if a plane section of a  surface $\Sigma \subset \mathbb{P}^3$ with a plane $H$ contains a node, this node is either a node of $\Sigma$ or a point of contact of $\Sigma$ and $H$. 
\begin{theorem}[Kummer]\label{kummer}
Let $\Sigma \subset \mathbb{P}^3$ be an irreducible quartic surface, not ruled in lines. Suppose that there exists a positive dimensional, algebraic family of planes $\mathcal{H}=\{H_t\}$ such that the section $H_t \cap \Sigma$ splits into a pair of conics. Then we are in one of the following cases.\\ \\
$(\mathbf{a})$ Suppose that the general plane of the family is not tangent to $\Sigma$; then we have several possibilities:
\begin{itemize} 
\item[($a1$)] two nodes are fixed for all the plane sections; in this case $\Sigma$ has a double conic and (at least) two isolated double points, and $\mathcal{H}$ coincides with the pencil of planes through these double points;
\item[($a2$)] three nodes are fixed for all the plane sections; in this case they belong to a line which is double for the surface, and $\mathcal{H}$ is the pencil of planes through this line;
\item[($a3$)] some of the double points of the sections always coincide; in this case either we are led to some special cases of $(a1)$, or we obtain a new type of surfaces, whose equation can be written in the form:
\begin{displaymath}  
\Phi^2=\alpha_1 \alpha_2 \alpha_3 \alpha_4,
\end{displaymath}
where $\Phi=0$ is a quadric and $\alpha_i=0$ are four coaxial planes. These surfaces contain two tacnodes, and $\mathcal{H}$ is the pencil of planes through these tacnodes.
\end{itemize}      
$(\mathbf{b})$ The general plane of $\mathcal{H}$ is tangent once to $\Sigma$; then we have the following possibilities:
\begin{itemize}
\item[($b1$)] none of the remaining nodes is fixed for all the plane sections of the family; in this case $\Sigma$ is a Steiner surface, and the family $\mathcal{H}$ is a family of tangent planes;
\item[($b2$)] one of the remaining nodes is fixed; in this case the surface $\Sigma$ contains a double conic and a node, and $\mathcal{H}$ is a family of tangent planes passing through this node. 
\end{itemize}
$(\mathbf{c})$ The general plane of $\mathcal{H}$ is bitangent to $\Sigma$; in this case $\Sigma$ contains a double conic. 
\end{theorem}
\begin{corollary} \label{kummer2}
Let $\Sigma \subset \mathbb{P}^3$ be an irreducible quartic surface which is not ruled in lines. Suppose that $\Sigma$ contains infinitely many conics; then one of the following holds:
\begin{itemize}   
\item[($i$)] $\Sigma$ is of the form 
\begin{displaymath}  
\Phi^2=\alpha_1 \alpha_2 \alpha_3 \alpha_4,
\end{displaymath}
where $\Phi=0$ is a quadric and $\alpha_i=0$ are four coaxial planes;
\item[($ii$)] $\Sigma$ has a double line;
\item[($iii$)] $\Sigma$ has a double conic;
\item[($iv$)] $\Sigma$ is a Steiner surface.
\end{itemize}
\end{corollary}
Let's now give a sketch of proof of Kummer's Theorem \ref{kummer}. \\ The sections of $\Sigma$ with the planes of $\mathcal{H}$ are quartic plane curves that split into two conics, hence they contains four nodes (possibly infinitely near). Let us consider first the case when no point is a point of contact. 
If none of the nodes is fixed, then $\Sigma$ contains a double curve of degree $4$; this implies that \emph{any} hyperplane section of $\Sigma$ has four nodes and hence is reducible, and this in turn implies that $\Sigma$ itself is reducible. If exactly one of the four nodes, say $P$, is fixed, than $\Sigma$ contains a double curve of degree $3$; hence all the plane sections through the point $P$ are reducible, and this implies either that $\Sigma$ itself is reducible or that it is a cone whose  vertex is $P$ (see [Ku, p.67] ). Then we can suppose that at least two nodes are fixed. If exactly two nodes are fixed, then the remaining two move and they form on $\Sigma$ a double conic; this is case $(a1)$. Suppose now that three nodes are fixed; therefore they necessarily lie on a line, that must be double for the surface: otherwise, the plane sections containing the nodes would be the union of a line and an irreducible plane cubic, instead of two conics; this is case $(a2)$. Finally, we must consider the case when the four nodes of the plane sections are not distinct; this happens in two cases: when we have some particular case of surface of type $(a1)$ (namely, when some of the isolated double points become infinitely near) and when the surface $\Sigma$ contains two tacnodes. In the latter case the equation of $\Sigma$ is $\Phi^2=\alpha_1 \alpha_2 \alpha_3 \alpha_4$, where $\Phi$ is a homogeneous polynomial of degree $2$ and the $\alpha_i$'s are four linear polynomial defining coaxial planes; the two tacnodes are given by the intersection of the axis of such planes with the quadric $\Phi=0$, and the pencil of planes through the tacnodes cuts out on $\Sigma$ a pencil of quartic curves that split into two bitangent conics. This is case $(a3)$. \\
Let us consider now the case when one of the nodes of the hyperplane sections is always a point of contact of the plane with the surface. If none of the remaining three nodes is fixed, then $\Sigma$ contains a double curve of degree $3$, and this curve can be given by a twisted cubic curve, by a conic and a line or by three lines concurrent at the same point; in the first two cases $\Sigma$ is ruled in lines, whereas in the last one we obtain a Steiner surface; it is indeed well known that the $\infty^2$ tangential plane sections of a Steiner surfaces are reducible into pairs of conics, and that conversely this property characterizes the Steiner surfaces (see \cite{Con39}); this is case $(b1)$. If one of the remaining nodes is fixed, then the other two form a double conic in $\Sigma$, and the algebraic system $\mathcal{H}$ is given by the the tangent planes to $\Sigma$ passing through the fixed node: this is case $(b2)$.\\
Finally, let us suppose that the general plane of $\mathcal{H}$ is bitangent to $\Sigma$; then none of the remaining nodes can be fixed, ant this in turn implies that $\Sigma$ contains a double conic. \\
The proof of Theorem \ref{kummer} in complete.\\ \\              
The following result is classically well-known.
\begin{proposition} \label{linearnorm} 
Let $\Sigma \subset \mathbb{P}^3$ be either a Steiner surface or an irreducible quartic surface with a double conic. Then $\Sigma$ is not linearly normal.
\end{proposition}  
\begin{proof}
The Steiner surface is the projection of the Veronese surface $V_4 \subset \mathbb{P}^5$ from a general line, hence it is not linearly normal. On the other hand, Segre showed in [Seg] (see also [Co39, Capitolo VI] e [Je16, Chapter III]) that any quartic surface with a double conic is the projection from an external point of a surface $\Phi_4 \subset \mathbb{P}^4$  which is complete intersection of two quadrics; hence it is not linearly normal.
\end{proof}
Now we are ready to show the following:
\begin{proposition} \label{non3}
The case $d=3$ never occurs.
\end{proposition}
\begin{proof}
Suppose that $d=3$ occurs. Then by Corollary \ref{infinite coniche} the bicanonical image $\Sigma$ is a surface in $\mathbb{P}^3$ which contains infinitely many conics, and that is not ruled in lines because of Proposition \ref{finite rette}. It follows that $\Sigma$ is one of the surfaces listed in Corollary \ref{kummer2}. Now, cases $(i)$ and $(ii)$ must be excluded because otherwise the conics would be cut out on $\Sigma$ by a \emph{pencil} of planes, contradicting Proposition \ref{nopencil}. On the other hand, cases $(iii)$ and $(iv)$ must be excluded because otherwise $\Sigma$ would be not linearly normal (Proposition \ref{linearnorm}), against the fact that the map $\phi$ is induced by the \emph{complete} bicanonical system $|2K|$.   
\end{proof}

\subsection{The case $d=4$} \label{d=4}

\begin{proposition}
\label{non in rette}
If $d=4$, then $\Sigma$ is not ruled in lines.
\end{proposition}
\begin{proof}
By contradiction, suppose that $\{L\}$ is a pencil of lines in $\Sigma$ and let $ \{\Lambda \}$ be its pullback in $S$. Therefore we have:
\begin{displaymath}
K\Lambda=\frac{1}{2}(2K)\Lambda=\frac{d}{2}=2.
\end{displaymath}
The Index Theorem gives $K^2 \Lambda^2 \leq (K \Lambda)^2$, that is 
$\Lambda^2 \leq 1$; since $K \Lambda=2, \; \Lambda^2$ is an
even integer, hence $\Lambda^2=0$. So $\{\Lambda\}$ is a genus $2$
pencil on $S$, therefore Corollary \ref{degbig} gives $d=2$ and this
is a contradiction.
\end{proof}
\begin{proposition} \label{non4}
If $X \in \mathfrak{B}$, then the case $d=4$ does not occur.
\end{proposition}
\begin{proof}
By contradiction, suppose $d=4$. Then the bicanonical image $\Sigma$ is a surface of degree $3$ in $\mathbb{P}^3$. If $\Sigma$ contains a double line, $\Sigma$ would be ruled in lines and this contradicts Proposition \ref{non in rette}; the same happens if $\Sigma$ contains a triple point. Therefore $\Sigma$ is a cubic surface with at worst rational double points, hence $K_{\Sigma}= -H$, where $H$ is a general hyperplane section of $\Sigma$. Now we apply the Hurwitz formula to the bicanonical morphism $\phi \colon S \longrightarrow \Sigma$; if $R$ is the ramification divisor, we obtain:
$$K \cong \phi^*K_{\Sigma}+R,$$
that is $R \cong 3K$. On the other hand, from the definition of $\mathfrak{B}$ it follows that the four curves $C_{\xi_i}, \; \; i= 0,1,2,3$ are smooth curves of genus $4$. The bicanonical divisors $2C_{\xi_i}$ correspond to hyperplane sections $H_i \subset \Sigma$. We claim that each $H_i$ is an irreducible, reduced plane cubic. 
Indeed, from the exact sequence 
\begin{equation} \label{biesatt}
0 \longrightarrow \mathcal{O}_S(C_{\xi_i}) \longrightarrow \mathcal{O}_S(2K_S) \longrightarrow \mathcal{O}_S(2K_S)|_{C_{\xi_i}} \longrightarrow 0
\end{equation}
it follows that $|2K|$ cuts out on $C_{\xi_i}$ a $g^2_6$ (which is complete if $i \neq 0$). Then the support of the divisor $H_i$ is an irreducible, nondegenerate curve in $\mathbb{P}^2$, and since $H_i$ is a hyperplane section and $\deg \Sigma=3$ our claim follows.
As $(2K)^2=12$, we obtain: 
\begin{displaymath}
\phi^*(H_{\xi_i})=2C_{\xi_i}, \quad 0 \leq i \leq 3,
\end{displaymath}
and this says that the four curves $\Co, \ldots, \Cthree$ are contained in the ramification divisor $R$ of $\phi$. Since $\Co + \cdots +\Cthree \cong 4K$, we obtain $h^0(S, \; R-4K)>0$. But this is impossible, because $R \cong 3K$.
\end{proof}

\subsection{The case $d=6$} \label{d=6}
\begin{proposition} \label{non6}
If $X \in \mathfrak{B}$, then the case $d=6$ does not occur.
\end{proposition}
\noindent We give two different proofs. \\ \\
\noindent \emph{First proof.}
Assume that $\phi$ has degree $6$. The bicanonical image $\Sigma$ is a linearly normal surface of degree $2$ in $\mathbb{P}^3$, hence a smooth quadric or a quadric cone.
But $\Sigma$ cannot be a smooth quadric; indeed, if $L$ is the pullback on $S$ of a ruling of $\Sigma$, we would have $L^2=0, \; KL=3$, and this is absurd because $KL+L^2$ must be an even integer. Hence $\Sigma$ is a quadric cone. We apply the Hurwitz formula to the bicanonical morphism $\phi \colon S \longrightarrow \Sigma$; if $R$ is the ramification divisor, we have:
\begin{displaymath}
K \cong \phi^*K_{\Sigma}+R,
\end{displaymath} 
that is $R \cong 5K$. Let $H_{\xi_i}$ be again the plane section of $\Sigma$ corresponding  to the bicanonical divisor $2C_{\xi_i}, \; \; i=0, \ldots, 3$. 
Since $|2K|$ cuts out on $C_{\xi_i}$ a $g^2_6$, it follows that $H_{\xi_i}$ is a smooth conic.
Therefore we have:
\begin{displaymath} 
\phi^*(H_{\xi_i})=2C_{\xi_i}, \quad 0 \leq i \leq 3,
\end{displaymath}
and this shows that the curves $C_{\xi_0}=K, \; C_{\xi_1}, \; C_{\xi_2}, \; C_{\xi_3}$ are contained in the ramification divisor $R$ with multiplicity $1$, thus  we can write:
\begin{displaymath}
R=\sum_{i=0}^3 C_{\xi_i}+D,
\end{displaymath}
where $D$ is an effective curve containing none of the curves $K, \; C_{\xi_1}, \; C_{\xi_2}, \; C_{\xi_3}$ in its support. But $R \cong 5K$, thus  we have:
\begin{displaymath}
D \cong 5K- \sum_{i=0}^3 C_{\xi_i} \cong K,
\end{displaymath}
hence $D=K$, a contradiction.
\begin{flushright}
$\square$
\end{flushright} 

\noindent \emph{Second Proof.}
We showed that $\Sigma$ is a quadric cone; projecting from its vertex, we obtain a rational map $S \dashrightarrow \mathbb{P}^1$, that is a linear pencil $|D|$ without fixed part such that $KD=3$ and $2K \cong 2D+G$, where $G$ is an effective divisor such that $KG=0$. This implies that every component of $G$ is a $(-2)-$curve, and since $X \in \mathfrak{B}$ we obtain $G=0$. Therefore $2D$ is linearly equivalent to $2K$, but $D$ is not linearly equivalent to $K$ because $p_g(S)=1$. This means that $D-K$ is a non trivial $2-$torsion element in $\textrm{Pic}(S)$, so we can consider the {\'e}tale double cover $ \pi \colon Y \longrightarrow S$ determined by $D-K$. Of course $Y$ is a minimal surface of general type and we have:
\begin{equation}
\begin{array}{ll}
K_Y^2=2K^2=6;\\
\chi(\mathcal{O}_Y)=2 \chi(\mathcal{O}_S)=2; \\
p_g(Y)=p_g(S)+h^0(S, \;D )=3,
\end{array}
\end{equation}
hence $q(Y)=2$. The Albanese fibration $\alpha \colon S \longrightarrow E$ of $S$ induces a fibration $f' \colon Y \longrightarrow E$.  Let $f \colon Y \longrightarrow B$ be the Stein factorization of $f'$, $g$ the genus of the general fibre of $f$ and $b$ the genus of $B$. Of course $b \geq 1$, because $B$ dominates $E$; on the other hand, the inequality $K_Y^2 \geq 8(g-1)(b-1)$ yields $b=1$, thus  $B$ is again an elliptic curve. Since $c_2(Y)= 18 \neq 0$, Zeuthen$-$Segre formula (see [BPV84, p.97]) implies that the fibration $f$ contains some singular fibres, hence $f$ is not locally trivial by the Grauert$-$Fischer theorem. Recall now that if $S$ is a surface and $f\colon S \longrightarrow B$ is a fibration over a smooth curve $B$ which is not locally trivial, it makes sense to consider the $\emph{slope}$ $\lambda(f)$ of $f$,  which is defined as the following ratio:
\begin{displaymath}
\lambda(f):=K^2_{S/B}/\textrm{deg}(f_*\omega_{S/B}).
\end{displaymath}
It is not difficult to see that $\lambda(f)$ is the unique number satisfying:
\begin{displaymath}
K_S^2=\lambda(f) \chi(\mO_S)+(8-\lambda(f))(b-1)(g-1).
\end{displaymath}
In order to conclude our second proof of Proposition \ref{non6}, we use the following result, whose proof can be found in [Xi87].
\begin{theorem}[Xiao] \label{slope}
Let $f \colon S \longrightarrow B$ be a minimal fibration, not locally trivial, with $g  \geq 2$. Then:
\begin{displaymath}
4(g-1)/g \leq \lambda(f) \leq 12
\end{displaymath}
and the second inequality becomes an equality if and only if every fiber of $f$ is smooth and reduced.\\
Moreover, if $\lambda(f) <4$, then $q(S)=g(B)$.
\end{theorem}
The slope of the fibration $f \colon Y \longrightarrow B$ is $\lambda(f)=K_Y^2 / \chi(\mathcal{O}_Y)=3$; hence the last statement in Theorem \ref{slope} gives $q(Y)=g(B)$. Since  $q(Y)=2$ and $g(B)=1$, we obtained a contradiction. \\
\begin{flushright}
$\square$
\end{flushright}
This concludes the proof of Theorem \ref{main2}.


\begin{thebibliography}{999999}

\bibitem[ADHPR93]{ADHPR93}
A. Aure, W. Decker, K. Hulek, S. Popescu, K. Ranestad: The geometry of bielliptic surfaces in $\mathbb{P}^4$, \emph{International Journal of Mathematics} $\boldsymbol{4}$, no. 6 (1993), 873-902. 

\bibitem[ACGH85]{ACGH85}
E. Arbarello, M. Cornalba, P. A. Griffiths and J. Harris:
\emph{Geometry of algebraic curves}, Springer-Verlag 1985.



\bibitem[At57]{At57}
M. F. Atiyah: Vector bundles over an elliptic curve, \emph{Proc. London Math. Soc.} $\boldsymbol{7}$ (1957), 414-452.




\bibitem[BPV84]{BPV84}
W. Barth, C. Peters, A. Van de Ven: \emph{Compact Complex
Surfaces}, Springer-Verlag 1984.

\bibitem[Be82]{Be82}
A. Beauville: L'inegalit{\'e} $p_g \geq 2q-4$ pour les surfaces de
type g{\'e}n{\'e}rale, \emph{Bull. Soc. Math. de France}
$\boldsymbol{110}$ (1982), 343-346.

\bibitem[Be88]{Be88}
A. Beauville: Annulation du $H^1$ et syst{\`e}mes paracanoniques sur le surfaces, 
\emph{J. reine angew. Math.} $\boldsymbol{388}(1988)$, 149-157. 

\bibitem[Be96]{Be96}
A. Beauville: \emph{Complex algebraic surfaces}, Cambridge
University Press 1996.


\bibitem[Bo73]{Bo73}
E. Bombieri: Canonical models of surfaces of general type,
\emph{Publ. IHES} $\boldsymbol{42}$ (1973), 171-219.

\bibitem[Bor03]{Bor03}
G. Borrelli: On the classification of surfaces of general type with non-birational bicanonical map and Du Val double planes, e-print AG/0312351

\bibitem[Ca81]{Ca81}
F. Catanese: On a class of surfaces of general type, in
\emph{Algebraic Surfaces}, CIME, Liguori (1981), 269-284.

\bibitem[Ca92]{Ca92} F. Catanese: Chow varieties, Hilbert Schemes and moduli space of surfaces of general type, \emph{Journal of Algebraic Geometry}, $\boldsymbol{1,4}$ (1992), 561-595.  


\bibitem[Ca99]{Ca99}
F. Catanese: Singular bidouble covers and the construction of interesting algebraic surfaces, \emph{Contemporary Mathematics} $\boldsymbol{241}$ (1999), 97-119.




\bibitem[CaCi91]{CaCi91}
F. Catanese and C. Ciliberto: Surfaces with $p_g=q=1$,
\emph{Symposia Math.} $\boldsymbol{32}$ (1991), 49-79.

\bibitem[CaCi93]{CaCi93}
F. Catanese and C. Ciliberto: Symmetric product of elliptic curves
and surfaces of general type with $p_g=q=1$, \emph{Journal of Algebraic
Geometry} $\boldsymbol{2}$ (1993), 389-411.

\bibitem[CCM98]{CCM98}
F. Catanese, C. Ciliberto and M. M. Lopes: Of the classification
of irregular surfaces of general type with non birational
bicanonical map, \emph{Trans. of the Amer. Math. Soc.}
$\boldsymbol{350}$ (1998), 275-308.

\bibitem[CD89]{CD89}
F. Catanese and O. Debarre: Surfaces with $K_S^2=2, p_g=1, q=0$,
\emph{J. Reine-Angew.-Math.} $\boldsymbol{395}$ (1989), 1-55.

\bibitem[Ch96]{Ch96} 
M. C. Chang: The number of components of Hilbert schemes, \emph{International Journal of Mathemathics} $\boldsymbol{7,3}$ (1996), 301-306.

\bibitem[CFM97]{CFM97}
C. Ciliberto, P. Francia, M. Mendes Lopes: Remarks on the bicanonical map of surfaces of general type, \emph{Math. Z.} $\boldsymbol{224}$ (1997) no.1, 137-166. 

\bibitem[CM02]{CM02}
C. Ciliberto, M. Mendes Lopes: On surfaces with $p_g=q=2$ and
non-birational bicanonical map, \emph{Adv. Geom.} $\boldsymbol{2}$ 
(2002), no. $3$, 281-300. 

\bibitem[CM02b]{CM02b}
C. Ciliberto, M. Mendes Lopes: On surfaces with $p_g=2, \; q=1$ and non-birational bicanonical map, \emph{Algebraic Geometry}, de Gruyter, Berlin (2002), 117-126.

\bibitem[Ci97]{Ci97}
C. Ciliberto: The bicanonical map for surfaces of general type,
\emph{Proc. of Symp. in Pure Mathematics} $\boldsymbol{62.1}$
(1997), 57-84.


\bibitem[Con39]{Con39}
F. Conforto: \emph{Le superficie razionali}, Zanichelli, Bologna
1939.

\bibitem[De82]{De82}
O. Debarre: In{\'e}galit{\'e}s num{\'e}riques pour les surfaces de type g{\'e}n{\'e}ral, \emph{Bull. Soc. Math. de France}
$\boldsymbol{110}$ (1982), 319-342.




\bibitem[En49]{En49}
F. Enriques: \emph{Le superficie algebriche}, Zanichelli, Bologna 1949.

\bibitem[Fr91]{Fr91}
P. Francia: On the base points of the bicanonical system, \emph{Symposia Math.} $\boldsymbol{32}$ (1991), 141-150.




\bibitem[Go31]{Go31}
L. Godeaux: Sur une surface alg{\'e}brique de genre zero and bigenre deux, \emph{Atti Acad. Naz. Lincei} $\boldsymbol{14}$ (1931), 479-481.   

\bibitem[Ha77]{Ha77}
R. Hartshorne: \emph{Algebraic Geometry}, Lecture Notes in Mathematics $\boldsymbol{52}$, Springer-Verlag 1977.

\bibitem[Ho77]{Ho77}
E. Horikawa: On algebraic surfaces with pencils of curves of genus
2, in Complex Analysis and Algebraic Geometry, a volume dedicated
to Kodaira, p.79-90, Cambridge 1977.



\bibitem[Je16]{Je16}
C. M. Jessop: \emph{Quartic surfaces with singular points},
Cambridge University Press, 1916.

\bibitem[Ku]{Ku}
E. Kummer:  {\"U}ber die Flachen vierten Grades auf Welchen
Schaaren von Kegelschnitte liegen, \emph{J. Reine-Angew.-Math.} $\boldsymbol{64}$ (1863), 66-76. 

\bibitem[Man97]{Man97}
M. Manetti: Iterated double covers and connected components of moduli spaces, \emph{Topology} $\boldsymbol{36,3}$ (1997), 745-764.

\bibitem[Mi84]{Mi84}
Y. Miyaoka: The maximum number of quotient singularities on
surfaces with given numerical invariants, \emph{Math. Ann.}
$\boldsymbol{268}$ (1984), 159-171.

\bibitem[MP01]{MP01}
M. Mendes Lopes, R. Pardini: The bicanonical map of surfaces with $p_g=0$ and $K^2 \geq 7$, \emph{Bull. London Math. Soc.} $\boldsymbol{33}$ no. 3 (2001), 265-274.  

\bibitem[MP03]{MP03}
M. Mendes Lopes, R. Pardini: The bicanonical map of surfaces with $p_g=0$ and $K^2 \geq 7$, II. \emph{Bull. London Math. Soc.} $\boldsymbol{35}$ no. 3 (2003), 337-343.  





\bibitem[Pol03]{Pol03}
F. Polizzi: Surfaces of general type with $p_g=q=1, \; K^2=8$ and bicanonical map of degree $2$, to appear in \emph{Transactions of the American Mathematical Society}.


\bibitem[Re88]{Re88}
I. Reider: Vector bundles of rank 2 and linear systems on
algebraic surfaces, \emph{Ann. of Math.} $\boldsymbol{127}$
(1988), 309-316.

\bibitem[Reid91]{Reid91}
M. Reid: Campedelli versus Godeaux, \emph{Symposia Math.} $\boldsymbol{32}$ (1991), 309-365.



\bibitem[Seg]{Seg}
C. Segre, Surfaces du quatri{\`e}me ordre {\`a} conique double, \emph{Math. Ann.} $\boldsymbol{24}$.

\bibitem[Xi85]{Xi85}
G. Xiao: Finitude de l' application bicanonique des surfaces de
type g{\'e}n{\'e}rale, \emph{Boll. Soc. Math. de France}
$\boldsymbol{113}$ (1985), 23-51.

\bibitem[Xi85b]{Xi85b}
G. Xiao: \emph{Surfaces fibr{\'e}es en courbes de genre deux}, Lecture
Notes in Mathematics $\boldsymbol{1137}$ (1985).

\bibitem[Xi87]{Xi87}
G. Xiao: Fibered algebraic surfaces with low slope,
\emph{Math. Ann.} $\boldsymbol{276}$ (1987), 449-466.

\bibitem[Xi90]{Xi90}
G. Xiao: Degree of the bicanonical map of a surface of general
type, \emph{Amer. J. of Math.} $\boldsymbol{112}$ $\boldsymbol{(5)}$, (1990) 309-316.




\end{thebibliography}
\end{document}